\RequirePackage{ifthen}
\newboolean{MPA}
\setboolean{MPA}{false}

\ifthenelse {\boolean{MPA}}
{
\documentclass[numbook]{svjour3} 
\smartqed
} {
\documentclass[11pt]{article}
\usepackage[margin=1in]{geometry}
}

\usepackage{graphicx}
\usepackage{amsmath}
\usepackage{amssymb}
\usepackage{hyperref}
\usepackage{color}

\usepackage{enumitem} 
\usepackage[capitalize,noabbrev]{cleveref}
\crefname{question}{Question}{Questions}
\crefname{step}{Step}{Steps}
\crefname{claim}{Claim}{Claims}
\crefname{problem}{Problem}{Problems}
\crefname{observation}{Observation}{Observations}

\usepackage[noadjust]{cite} 

\usepackage{mathrsfs} 

\usepackage{xspace} 

\usepackage{soul} 

\usepackage[makeroom]{cancel} 

\usepackage{thmtools}
\usepackage{thm-restate}

\newcommand{\R}{\mathbb R}
\newcommand{\Z}{\mathbb Z}
\newcommand{\Q}{\mathbb Q}

\newcommand{\badone}{x^1}
\newcommand{\badtwo}{x^2}
\newcommand{\badoney}{y^1}
\newcommand{\badtwoy}{y^2}

\let\st\relax
\DeclareMathOperator{\st}{s.t.}

\DeclareMathOperator{\spn}{span}
\DeclareMathOperator{\size}{size}
\DeclareMathOperator{\width}{width}

\newcommand{\genobj}{f}

\renewcommand{\P}{\mathcal P}
\newcommand{\B}{\mathcal B}

\renewcommand{\L}{\mathcal L}

\renewcommand{\S}{\mathcal S}

\newcommand{\E}{\mathcal E}

\renewcommand{\H}{\mathcal H}

\newcommand{\direction}{d}

\newcommand{\floor}[1]{\lfloor#1\rfloor}

\newcommand{\ceil}[1]{\lceil#1\rceil}

\newcommand{\abs}[1]{\lvert#1\rvert}

\newcommand{\norm}[1]{\lVert#1\rVert_2}

\newcommand{\pare}[1]{\left(#1\right)}
\newcommand{\bra}[1]{\left\{#1\right\}}
\newcommand{\sbra}[1]{\left[#1\right]}

\newcommand{\transp}{\mathsf T}

\newcommand{\constlen}{2^{p(p-1)/4}}

\usepackage{comment}
\specialcomment{knownproof}{\color{blue}}{\color{black}}
\excludecomment{knownproof}

\newcounter{step}

\ifthenelse {\boolean{MPA}}
{

\newenvironment{prf}[1][]
{\begin{proof}}
{\qed \end{proof}}

\newenvironment{prfc}[1][]
{\begin{proof}[#1]}
{\qed \end{proof}}

\newenvironment{prfh}[1][]
{\begin{proof}}
{\end{proof}}

\newcounter{claim} 
\renewenvironment{claim}[1][]
{\refstepcounter{claim} \begin{trivlist} \item[] {\bf Claim~\theclaim}\space#1 \itshape}
{\end{trivlist}}

\newenvironment{cpf}
{\begin{trivlist} \item[] {\em Proof of claim }}
{$\hfill\diamond$ \end{trivlist}}

\journalname{Mathematical Programming A}

} {

\usepackage{amsthm}
\newtheorem{theorem}{Theorem}
\newtheorem{proposition}{Proposition}
\newtheorem{lemma}{Lemma}

\newenvironment{prf}[1][]
{\begin{proof}}
{\end{proof}}

}

\newtheorem{observation}{Observation}

\begin{document}


\title{The Mixed Integer Trust Region Problem}

\ifthenelse {\boolean{MPA}}
{
\titlerunning{An Approximation Algorithm for Indefinite MIQP}

\author{Alberto Del Pia}
\institute{Alberto~Del~Pia \at
              Department of Industrial and Systems Engineering 
              \& Wisconsin Institute for Discovery \\
              University of Wisconsin-Madison, Madison, WI, USA \\
              \email{delpia@wisc.edu}}
}
{
\author{Alberto Del Pia
\thanks{Department of Industrial and Systems Engineering \& Wisconsin Institute for Discovery,
             University of Wisconsin-Madison, Madison, WI, USA.
             E-mail: {\tt delpia@wisc.edu}.}}
}

\date{September 17, 2024}

\pagenumbering{gobble} 

\maketitle

\pagenumbering{arabic}

\begin{abstract}
In this paper we consider the problem of minimizing a general quadratic function over the mixed integer points in an ellipsoid.
This problem is strongly NP-hard, NP-hard to approximate within a constant factor, and optimal solutions can be irrational.
In our main result we show that an arbitrarily good solution can be found in polynomial time, if we fix the number of integer variables.
This algorithm provides a natural extension to the mixed integer setting, of the polynomial solvability of the trust region problem proven by Ye, Karmarkar, Vavasis, and Zippel.
As a result, our findings pave the way for designing efficient trust region methods for mixed integer nonlinear optimization problems.
The techniques that we introduce are of independent interest and can be used in other mixed integer nonlinear optimization problems.
As an example, we consider the problem of minimizing a general quadratic function over the mixed integer points in a polyhedron.
For this problem, we show that a solution satisfying weak bounds with respect to optimality can be computed in polynomial time, provided that the number of integer variables is fixed.
It is well known that finding a solution satisfying stronger bounds cannot be done in polynomial time, unless P=NP.
\ifthenelse {\boolean{MPA}}
{
\keywords{trust region problem \and mixed integer quadratic programming \and approximation algorithm \and polynomial time algorithm \and ellipsoid constraint}
} {}
\end{abstract}

\ifthenelse {\boolean{MPA}}
{}{
\emph{Key words:} trust region problem; mixed integer quadratic programming; approximation algorithm; polynomial time algorithm; ellipsoid constraint
}



\section{Introduction}
\label{sec intro}

\emph{Trust region methods} form an important class of algorithms for nonlinear optimization problems.
The appeal of these methods lies in their general applicability, the elegant theory, and their practical success.
We refer the reader to the book \cite{ConGouToi00book} for a thorough introduction.
In trust region methods, a quadratic function is used to model the nonlinear objective function, and this model is expected to be accurate (the model is ``trusted'') only in a neighborhood of the current iterate $c$.
To obtain the next iterate, first an affine change of variables is preformed to map $c$ to the origin, and then the model quadratic function is minimized over a neighborhood of the origin.
If the 2-norm is used to define the neighborhood, the problem takes the form
\begin{align}
\label[problem]{prob TR}
\tag{TR}
\begin{split}
\min & \quad x^\transp H x + h^\transp x \\
\st & \quad x^\transp x \le 1.
\end{split}
\end{align}
Here $H$ is a symmetric matrix in $\Q^{n \times n}$ and $h \in \Q^n$.
\cref{prob TR} 
is often referred to, in the literature, as the \emph{trust region} (sub)problem.

One key reason behind the success of trust region methods is the fact that \cref{prob TR} can be solved efficiently.
In fact, while optimal solutions
of \cref{prob TR} can be irrational (see \cref{obs TR irrational}),
an arbitrarily good solution can be found in polynomial time.
This celebrated result is due to Ye \cite{Ye92b}, Karmarkar \cite{Kar89}, and Vavasis and Zippel \cite{VavZip90}, and a formal statement is given in \cref{lem VavZip90}.
Tractability results for extensions of \cref{prob TR} that include additional linear or quadratic constraints can be found, for example, in \cite{FuLuoYe98,YeZha03,BieMic14,Bie16,ConLoc17}.

Trust region methods have been introduced also for mixed integer nonlinear optimization problems
(see, e.g., \cite{ExlSch07,ExlerThesis13,NewbyThesis13,NewAli15}).
For this more general type of optimization problems, 
if the 2-norm is used to define the neighborhood, 
the next iterate can be obtained by solving the \emph{mixed integer trust region} (sub)problem
\begin{align}
\label[problem]{prob MITR}
\tag{MITR}
\begin{split}
\min & \quad x^\transp H x + h^\transp x \\
\st & \quad x^\transp x \le 1 \\
& \quad x \in \Pi_p(b^1,\dots,b^n) + \{c\}.
\end{split}
\end{align}
Here $H$ is a symmetric matrix in $\Q^{n \times n}$, $h,c \in \Q^n$, $b^1,\dots,b^n \in \Q^n$ are linearly independent, $p \in \{0,1,\dots,n\}$, and $\Pi_p(b^1, \dots , b^n)$ denotes the \emph{mixed integer lattice} 
$$
\Pi_p(b^1, \dots , b^n) := \bra{\sum_{i=1}^n \mu_i b^i : \mu_i \in \Z \ \forall i=1,\dots,p, \ \mu_i \in \R \ \forall i=p+1,\dots,n},
$$
which is the image of the set $\Z^p \times \R^{n-p}$ under an affine change of variables.
As opposed to \cref{prob TR}, we will see in \cref{sec hardness} that \cref{prob MITR} is strongly NP-hard and NP-hard to approximate within a constant factor, even in the pure integer case $p=n$.
Even deciding whether the feasible region of \cref{prob MITR} is nonempty, is NP-hard.

The relationship between the computational complexity of \cref{prob TR,prob MITR} resembles the situation in linear optimization, where linear programs can be solved in polynomial time, while mixed integer linear programs are strongly NP-hard and NP-hard to approximate within a constant factor.
The key link
between these two results was proven by Lenstra \cite{Len83}: mixed integer linear programs can be solved in polynomial time, if the number of integer variables is fixed (but the number of continuous variables does not need to be fixed).
In our main result
we provide the same type of 
link
between \cref{prob TR,prob MITR}: we prove that an arbitrarily good solution for \cref{prob MITR} can be found in polynomial time, if the number $p$ of integer variables is fixed.
The key difference between Lenstra's result and ours, besides the different types of optimization problems considered, is that here we deal with approximate solutions rather than global optimal solutions, but this is unavoidable due to the fact that optimal solutions to \cref{prob TR,prob MITR} are generally irrational.
In particular, our findings pave the way for designing efficient trust region methods for mixed integer nonlinear optimization problems.


\medskip
\noindent
\textbf{Ellipsoid-constrained mixed integer quadratic programming.}
Our main result actually holds for a problem that is more general than \cref{prob MITR}. We call it the \emph{ellipsoid-constrained mixed integer quadratic programming} problem, and we define it as follows:
\begin{align}
\label[problem]{prob E-MIQP}
\tag{E-MIQP}
\begin{split}
\min & \quad x^\transp H x + h^\transp x \\
\st & \quad (x-c)^\transp Q (x-c) \leq 1 \\
& \quad x \in \Z^p \times \R^{n-p}.
\end{split}
\end{align}
Here $H$ is a symmetric matrix in $\Q^{n \times n}$, $Q$ is a symmetric positive definite matrix in $\Q^{n \times n}$, $h,c \in \Q^n$, and $p \in \{0,\dots,n\}$.
There are two key differences with \cref{prob MITR}: the ball constraint has been replaced with a more general ellipsoid constraint, and the mixed integer lattice $\Pi_p(b^1,\dots,b^n) + \{c\}$ has been replaced with the ``cleaner'' mixed integer lattice $\Z^p \times \R^{n-p}$.

As we will see in \cref{sec relationship}, \cref{prob MITR} is essentially the special case of \cref{prob E-MIQP} where the matrix $Q$ is of the form $Q = B^\transp B$ for some invertible matrix $B \in \Q^{n \times n}$.
On the other hand, \cref{prob E-MIQP} is indeed more general than \cref{prob MITR}, in the sense that there are instances of \cref{prob E-MIQP} for which there exists no rational affine transformation that maps them to \cref{prob MITR}.

\medskip
\noindent
\textbf{An approximation algorithm for \cref{prob E-MIQP}.}
To state our main results, we first define 
the concepts of $\epsilon$-approximate solutions and of size.
Consider 
a general optimization problem $\inf \bra{ \genobj(x) : x \in \S}.$
Let $\genobj_{\inf}$ and $\genobj_{\sup}$ denote the infimum and the supremum of~$\genobj(x)$ on 
$\S$.
For $\epsilon \in [0,1]$, we say that $x^\diamond \in \S$ is an \emph{$\epsilon$-approximate solution} if
\begin{equation*}
\genobj(x^\diamond) - \genobj_{\inf} \le \epsilon \cdot (\genobj_{\sup} - \genobj_{\inf}).
\end{equation*}
Note that only optimal solutions are $0$-approximate solutions, while any point in $\S$ is a $1$-approximate solution.
Our definition of approximation 
has
some well-known 
invariance properties which make it a natural choice for the problems considered in this paper.
For instance, it is preserved under dilation and translation of the objective function, and it is insensitive to affine transformations of the objective function and of the feasible region, like for example changes of basis.
This definition has been used in our earlier work, and in numerous other places such as~\cite{NemYud83,Vav92c,BelRog95,KleLauPar06}.
The \emph{size} (also known as \emph{bit size}, or \emph{length}) of rational numbers, vectors, matrices, constraints, and optimization problems, denoted by $\size(\cdot)$, is the standard one in mathematical programming (see, e.g., \cite{SchBookIP,ConCorZamBook}), and is essentially the number of bits required to encode such objects.
We are now ready to state the first main result of this paper.

\begin{theorem}[Approximation algorithm for \cref{prob E-MIQP}]
\label{th MITR algorithm}
There is an algorithm which, for any $\epsilon \in (0,1]$, either finds an $\epsilon$-approximate solution to \cref{prob E-MIQP}, or correctly detects that the feasible region is empty.
The running time of the algorithm is polynomial in the size of \cref{prob E-MIQP} and in $1/\epsilon$, provided that the number $p$ of integer variables is fixed.
\end{theorem}



\cref{th MITR algorithm} implies that it is significantly simpler to optimize a quadratic function over mixed integer points in an ellipsoid rather than over mixed integer points in a polyhedron, a behavior previously only known in the pure continuous setting \cite{Vav93,ParVav91}.
In turn, this suggests that it might be advantageous to design trust region methods for mixed integer nonlinear optimization problems that are based on the 2-norm, or ellipsoidal norms, rather than polyhedral norms.

Even though \cref{prob MITR} is a special case of \cref{prob E-MIQP}, it still plays a central role in our path to prove \cref{th MITR algorithm}.
In \cref{sec TR,sec MITR}, we study \cref{prob TR,prob MITR}, and we prove a number of algorithmic and structural results of independent interest.
The end goal of these two sections is to prove \cref{prop MITR dist points}, which introduces an algorithm for \cref{prob MITR} that either finds two feasible solutions whose difference in objective value is at least $\max \bra{\norm{H},\norm{h}}$, for $\epsilon$ that goes to zero, or finds a nonzero vector $\direction$ along which the feasible region is flat.
In \cref{sec E-MIQP alt}, we get back to the more general \cref{prob E-MIQP}.
Mapping the ellipsoid defining its feasible region to a new ellipsoid that is arbitrarily close to the unit ball, allows us to employ \cref{prop MITR dist points} to obtain a flatness result for \cref{prob E-MIQP}.
This result, presented in \cref{prop E-MIQP flat}, gives an algorithm that either finds an $\epsilon$-approximate solution, or finds a nonzero vector $d$ along which the feasible region is flat.
In the second case, we are able to map the intersection of the feasible region with a hyperplane, to the feasible region of a new \cref{prob E-MIQP} with one fewer integer variable, and as a result we complete our proof of \cref{th MITR algorithm}.
As we will see, this last step is where the additional level of generality of \cref{prob E-MIQP} over \cref{prob MITR} plays a key role.
For more details on the ideas behind our proofs and techniques used, we refer the reader to the subsequent sections, where some proof overviews are provided to complement the complete and detailed proofs of the main results.


\medskip
\noindent
\textbf{A weak approximation algorithm for \cref{prob MIQP}.}
The techniques that we introduce to prove \cref{th MITR algorithm} are of independent interest and can be used in other mixed integer nonlinear optimization problems. 
To demonstrate this claim, we consider the \emph{mixed integer quadratic programming} problem
\begin{align}
\label[problem]{prob MIQP}
\tag{MIQP}
\begin{split}
\min & \quad x^\transp H x + h^\transp x \\
\st & \quad Wx \le w \\
& \quad x \in \Z^p \times \R^{n-p}.
\end{split}
\end{align}
Here $H$ is a symmetric matrix in $\Q^{n \times n}$, $h \in \Q^n$, $W \in \Q^{m \times n}$, $w \in \Q^m$, and $p \in \{0,\dots,n\}$.
\cref{prob MIQP} is NP-hard even in the pure continuous setting \cite{ParVav91}.
In the general mixed integer case, it is NP-hard to decide whether the feasible region of \cref{prob MIQP} is nonempty.
However, this feasibility problem can be solved in polynomial time if the number $p$ of integer variables is fixed \cite{Len83}, thus we now focus on this case and assume $p$ is fixed in the subsequent discussion.
In \cite{dP23bMPA}, the author gives an algorithm that finds an $\epsilon$-approximate solution to \cref{prob MIQP} in time polynomial in the size of \cref{prob MIQP} and in $1/\epsilon$, provided that the rank of $H$ is fixed.
In the general case, where 
the rank of $H$ is not fixed, 
no such algorithm can exist, unless P=NP \cite{Vav92c}.
In fact, unless P=NP, there is a constant $\epsilon \in (0, 1)$ such that no algorithm can find an $\epsilon$-approximate solution to \cref{prob MIQP} in polynomial time \cite{BelRog95}.
As a result, we cannot hope to approximate \cref{prob MIQP}
in polynomial time unless we are willing to accept an approximation factor that tends to one asymptotically as the problem gets larger \cite{Vav93}.
In \cref{sec MIQP}, we show how the algorithmic and structural results on \cref{prob MITR}, and in particular \cref{prop MITR dist points}, can be used to design a polynomial time approximation algorithm for \cref{prob MIQP} that satisfies a weak approximation bound of this sort.
\begin{theorem}[Weak approximation algorithm for \cref{prob MIQP}]
\label{th MIQP algorithm}
Consider \cref{prob MIQP} and assume there exists $\bar f \in \R$ such that every feasible solution has objective value at least $\bar f$.
There is an algorithm which either finds an $\epsilon$-approximate solution to the problem,
where the approximation factor is $\epsilon = 1-\Theta(n^{-2})$, or correctly detects that the feasible region is empty.
The running time of the algorithm is polynomial in the size of \cref{prob MIQP}, provided that the number $p$ of integer variables is fixed.
\end{theorem}

\cref{th MIQP algorithm} provides an extension to the mixed integer setting of a result by Vavasis in the pure continuous setting \cite{Vav93}. 
Note that the assumption in \cref{th MIQP algorithm} that the objective function is lower bounded on the feasible region cannot be relaxed, unless P=NP \cite{dP23bMPA}.
%

We remark that 
the main emphasis of this paper lies on the theoretical computational complexity of the algorithms presented. 
We do not refer to practically efficient implementations of the algorithms, which we believe should be studied in the future.
Before proceeding with the proofs of our main results, in \cref{sec relationship} we discuss the relationship between \cref{prob MITR,prob E-MIQP}, and in \cref{sec hardness} we study the hardness of the Problems~\ref{prob TR}, \ref{prob MITR}, and \ref{prob E-MIQP}.

\section{Relationship between \cref{prob MITR,prob E-MIQP}}
\label{sec relationship}

The next result shows that \cref{prob MITR} is essentially the special case of \cref{prob E-MIQP} where the matrix $Q$ is of the form $Q = B^\transp B$ for some invertible matrix $B \in \Q^{n \times n}$.

\begin{observation}
\label{obs special case}
Up to a rational affine transformation, \cref{prob MITR} is the special case of \cref{prob E-MIQP} where the matrix $Q$ is of the form $Q = B^\transp B$ for some invertible matrix $B \in \Q^{n \times n}$.
\end{observation}

\begin{prf}
Consider \cref{prob MITR} and perform the change of variables $y := B^{-1} (x - c)$, where $B$ is the invertible matrix with columns $b^1,\dots,b^n$.
The optimization problem that we obtain can be written in the form
\begin{align*}
\gamma + \min & \quad y^\transp H' y + {h'}^\transp y \\
\st & \quad (y-c')^\transp Q (y-c') \leq 1 \\
& \quad y \in \Z^p \times \R^{n-p},
\end{align*}
where $\gamma:= c^\transp H c + h^\transp c$, $H':= B^{\transp} H B$, $h':= 2 B^\transp H^\transp c + B^\transp h$, $Q:= B^\transp B$, and $c':= -B^{-1}c$.

Next, consider \cref{prob E-MIQP} and assume that the matrix $Q$ is of the form $Q = B^\transp B$ for some invertible matrix $B \in \Q^{n \times n}$.
We perform the change of variables $y := B (x - c)$ and obtain the optimization problem
\begin{align*}
\gamma + \min & \quad y^\transp H' y + {h'}^\transp y \\
\st & \quad y^\transp y \le 1 \\
& \quad y \in \Pi_p(b^1,\dots,b^n) + \{c'\},
\end{align*}
where $\gamma:= c^\transp H c + {h}^\transp c$, $H' := B^{-\transp} H B^{-1}$, 
$h' := 2 B^{-\transp} H^\transp c + B^{-\transp} h$,
$b^1,\dots,b^n$ are the columns of $B$,
and $c' := -Bc$.
\end{prf}

Next, we show that \cref{prob E-MIQP} is indeed more general than \cref{prob MITR}, by presenting instances of \cref{prob E-MIQP} for which there exists no rational affine transformation that maps them to \cref{prob MITR}.
Clearly, there is no rational affine transformation that maps the
$1$-dimensional ellipsoid $\bra{x \in \R : 2 x_1^2 \leq 1}$ to the $1$-dimensional ball $\bra{x \in \R : x_1^2 \leq 1}$.
However, this simple example could be written as \cref{prob MITR}, if the constraint $x^\transp x \le 1$ is replaced with $x^\transp x \le r$, for some positive $r \in \Q$.
Next, we give a more interesting $2$-dimensional example which shows that \cref{prob E-MIQP} is also more general than the extension of \cref{prob MITR} with a constraint of the form $x^\transp x \le r$, for some positive $r \in \Q$.

\begin{observation}
There exists no rational affine transformation that maps the $2$-dimensional ellipsoid $\bra{x \in \R^2 : 2 x_1^2 + x_2^2 \leq 1}$ to a $2$-dimensional ball $\bra{x \in \R^2 : x_1^2 + x_2^2 \leq r}$, for some positive $r \in \Q$.
\end{observation}

\begin{prf}
Assume, for a contradiction, that there exists a rational affine transformation as in the statement.
The affine transformation is of the form $y = B(x-c)$, for some invertible matrix $B \in \Q^{2 \times 2}$, and $c \in \Q^2$.
We must have $c=0$, since the centers of the ellipsoid and of the ball are both the origin.
We then have
\begin{align*}
& x^\transp
\begin{pmatrix}
2r & 0 \\
0 & r
\end{pmatrix}
x
=
y^\transp
y 
=
x^\transp B^\transp
Bx.
\end{align*}
If we denote by $b_{ij}$ the element of $B$ in position $(i,j)$, these elements must satisfy the system
\begin{align}
\label{eq system}
\begin{split}
& b_{11}^2 + b_{21}^2 = 2r \\
& b_{12}^2 + b_{22}^2 = r \\
& b_{11} b_{12} + b_{21} b_{22} = 0.
\end{split}
\end{align}
From the first two equations in \eqref{eq system}, we derive 
\begin{align}
\label{eq 2}
2 = \frac{b_{11}^2 + b_{21}^2}{b_{12}^2 + b_{22}^2}.
\end{align}
From the third equation in \eqref{eq system}, we obtain
\begin{align*}
b_{11}^2 b_{12}^2 = b_{21}^2 b_{22}^2.
\end{align*}
Since $B$ is invertible, $b_{12}$ and $b_{22}$ cannot both equal zero.
If $b_{12} \neq 0$, we multiply the right hand of \eqref{eq 2} by $b_{12}^2 / b_{12}^2$ and conclude
\begin{align*}
2 = \frac{b_{11}^2 b_{12}^2 + b_{21}^2 b_{12}^2}{b_{12}^2 (b_{12}^2 + b_{22}^2)}
= \frac{b_{21}^2 b_{22}^2 + b_{21}^2 b_{12}^2}{b_{12}^2 (b_{12}^2 + b_{22}^2)}
= \frac{b_{21}^2 (b_{12}^2 + b_{22}^2)}{b_{12}^2 (b_{12}^2 + b_{22}^2)}
= \frac{b_{21}^2}{b_{12}^2}.
\end{align*}
Hence, $\sqrt{2}$ is the rational number $\abs{b_{21}/b_{12}}$, a contradiction.
On the other hand, if $b_{22} \neq 0$, we multiply the right hand of \eqref{eq 2} by $b_{22}^2 / b_{22}^2$ and conclude
\begin{align*}
2 = \frac{b_{11}^2 b_{22}^2 + b_{21}^2 b_{22}^2}{b_{22}^2 (b_{12}^2 + b_{22}^2)}
= \frac{b_{11}^2 b_{22}^2 + b_{11}^2 b_{12}^2}{b_{22}^2 (b_{12}^2 + b_{22}^2)}
= \frac{b_{11}^2 (b_{12}^2 + b_{22}^2)}{b_{22}^2 (b_{12}^2 + b_{22}^2)}
= \frac{b_{11}^2}{b_{22}^2}.
\end{align*}
Hence, $\sqrt{2}$ is the rational number $\abs{b_{11}/b_{22}}$, a contradiction.
\end{prf}

\section{Hardness of mixed integer trust region problems}
\label{sec hardness}

In this section, we present some hardness results for Problems~\ref{prob TR}, \ref{prob MITR}, and \ref{prob E-MIQP} that we mentioned in the introduction.
We recall that these three problems are listed in increasing order of generality.
We start with a well-known fact about \cref{prob TR}.

\begin{observation}
\label{obs TR irrational}
There are instances of \cref{prob TR} where the optimal value is irrational and no optimal solution is rational, even if $n=2$.
\end{observation}

\begin{prf}
Consider the instance of \cref{prob TR} with $H:=\begin{pmatrix} 2 & 1 \\ 1 & 0 \end{pmatrix}$, and $h:=0$, i.e., $\min\bra{2x_1^2 + 2x_1x_2 : x^\transp x \le 1}$.
It can be checked that the smallest eigenvalue of $H$ is $\lambda_{\min} := 1-\sqrt 2$.
It is well known that $\lambda_{\min} = \min\bra{x^\transp H x : x^\transp x = 1}$ and, since $\lambda_{\min} \le 0$, we also have $\lambda_{\min} = \min \bra{ x^\transp H x : x^\transp x \le 1}$.
Hence the optimal value of this instance is $\lambda_{\min}$.
This instance has a unique optimal solution, which is the unit length eigenvector $(1-\sqrt 2,1)/\norm{(1-\sqrt 2,1)}$ associated with the eigenvalue $\lambda_{\min}$.
\end{prf}

\cref{obs TR irrational} also holds for the more general Problems~\ref{prob MITR}, and \ref{prob E-MIQP}, even if $p=0$ and $n=2$.
These results imply that no polynomial time algorithm, which stores numbers in binary encoding, can find an optimal solution to these problems.
Next, we discuss the hardness of \cref{prob MITR}.
To do so, we need a lemma on the approximation of square roots.
While the square root of a rational number can be irrational, it is well known that it can be approximated efficiently via binary search.
We formalize this simple result in the next lemma.
The proof is given for completeness.

\begin{lemma}
\label{lem sqrt approx}
Let $\alpha \in \Q$ be positive.
There is an algorithm which, for any positive integer $k$, finds nonnegative $l,u \in \Q$ such that $l \le \sqrt{\alpha} \le u$ and $u-l = \max \bra{1,\alpha} / 2^k$.
The running time of the algorithm is polynomial in $\size(\alpha), k$.
\end{lemma}

\begin{prf}
Let $l_0:=0$, $u_0:=\max \bra{1,\alpha}$.
Our algorithm consists of $k$ iterations $i=1,\dots,k$.
At iteration $i$, set $m_i := (l_{i-1}+u_{i-1})/2$.
Compare $m_i^2$ with $\alpha$.
If $m_i^2 \le \alpha$, set $l_{i}:=m_i$ and $u_{i}:=u_{i-1}$.
If $m_i^2 > \alpha$, set $l_{i}:=l_{i-1}$ and $u_{i}:=m_i$.
When iteration $k$ is completed, 
return $l:=l_k$ and $u:=u_k$.

To show correctness of the algorithm, note that $\sqrt{\alpha} \in [0, \max \bra{1,\alpha}]$.
The algorithm performs binary search on this interval, comparing, at each iteration, the midpoint $m_i$ of the previous interval with $\sqrt{\alpha}$, and updating the interval accordingly.
In fact, we have $m_i \lesseqgtr \sqrt{\alpha}$ if and only if $m_i^2 \lesseqgtr \alpha$.
Hence, for every $i=0,1,\dots,k$, we have $\sqrt{\alpha} \in [l_i,u_i]$.
Furthermore, for every $i=0,1,\dots,k$, we have $u_i-l_i = \max \bra{1,\alpha} / 2^i$.
\end{prf}

\begin{theorem}
\label{th MITR hardness}
\cref{prob MITR} is strongly NP-hard and it is NP-hard to find an $\epsilon$-approximate solution with $\epsilon <1/17$, even if $p=n$ and the problem is feasible.
\end{theorem}

\begin{prf}
The Max-Cut problem is strongly NP-hard \cite{GarJohSto76} and it is NP-hard to find an $\epsilon$-approximate solution, with $\epsilon >16/17$ \cite{TreSorSudWil00}.
It is well known that the Max-Cut problem on a graph $G = (V,E)$ with $V = \{1,2,\dots,n\}$ and nonnegative weights $w_{uv} \in \Z$, for $uv \in E$, can be formulated as
\begin{align*}
\max & \quad \sum_{uv \in E} \frac{w_{uv}}{2} (1-x_u x_v) \\
\st & \quad x \in \{-1,1\}^n.
\end{align*}
Observe that $\{-1,1\}^n = \{x \in 2\Z^n + \{1\} : x^\transp x \le r\},$ for any $r \in [n,n+8)$.
We use \cref{lem sqrt approx} with $\alpha:=n$ and $k := n$.
We obtain $l,u \in \Q$ such that $l \le \sqrt{n} \le u$ and $u-l = n / 2^n$.
Therefore, $u-\sqrt{n} \le n / 2^n < \sqrt{n+8} - \sqrt{n}$, where the last inequality can be checked with little effort.
Hence we have $u \in [\sqrt{n},\sqrt{n+8})$, thus $u^2 \in [n,n+8)$.
Applying the change of variables $y := x / u$ to the above formulation of the Max-Cut problem, we obtain the following special case of \cref{prob MITR}:
\begin{align*}
\max & \quad \sum_{uv \in E} \frac{w_{uv}}{2} (1-u^2 y_u y_v) \\
\st & \quad y^\transp y \le 1 \\
& \quad y \in \frac{2}{u}\Z^n + \bra{\frac{1}{u}}.
\end{align*}
\end{prf}


%

Our next goal is to show that it is NP-hard to decide whether the feasible region of \cref{prob MITR} is nonempty.
To do so, we first define the concept of lattice.
Given linearly independent vectors $b^1, \dots , b^p \in \R^n$, we define the \emph{lattice} 
$$
\Lambda(b^1, \dots , b^p)
 := \bra{\sum_{i=1}^p \mu_i b^i : \mu_i \in \Z \ \forall i=1,\dots,p}.
$$
Then, we consider the \emph{decisional version of the closest vector problem (CVP):} 
Given linearly independent vectors $b^1, \dots , b^p \in \Q^n$, 
a target vector $t \in \Q^n$, and $r \in \Q$, determine whether there exists $y \in \Lambda(b^1, \dots , b^p)$ with $\norm{y-t} \le r$.
This problem is well known to be NP-hard (see theorem 3.1 in \cite{MicGolBook}).
We first observe that the problem is NP-hard, even if $p=n$.

\begin{lemma}
\label{lem full rank CVP}
The decisional version of the CVP is NP-hard, even if $p=n$.
\end{lemma}

\begin{prf}
We give a polynomial reduction from the decisional version of the CVP.
For $j=p+1,\dots,n$, we can find in polynomial time a vector $b^j \in \Q^n$ orthogonal to $b^1,\dots,b^{j-1}$ and with $\norm{b^{j}} > \norm{t} + r$.
This can be done by appropriately rescaling a nonzero solution to the system of linear equations ${b^i}^\transp x = 0$, for $i=1,\dots,j-1$.
Clearly $b^1, \dots , b^n$ are linearly independent. 
It is simple to check that the decisional version of the CVP with data $b^1, \dots , b^p, t, r$ is equivalent to the decisional version of the CVP with data $b^1, \dots , b^n, t, r$.
\end{prf}

\begin{theorem}
\label{th MITR fesibility hardness}
It is NP-hard to decide whether the feasible region of \cref{prob MITR} 
is nonempty, even if $p=n$.
\end{theorem}


\begin{prf}
We give a polynomial reduction from the decisional version of the CVP with $p=n$.
Let 
$\S := \{y \in \Lambda(b^1, \dots , b^n) : \norm{y-t} \le r\}$.
Define $x := (y-t)/r$.
We obtain that $\S$ is nonempty if and only if the following set is nonempty:
\begin{align*}
\bra{x \in \Lambda\pare{\frac{b^1}{r}, \dots , \frac{b^n}{r}} : \norm{x} \le 1} - \bra{\frac{t}{r}}.
\end{align*}
The latter set is the feasible region of \cref{prob MITR} with
$b^i := b^i/r$ for $i=1,\dots,n$, $c := -t/r$,
and $p:=n$.
\end{prf}

Since \cref{prob MITR} is a special case of \cref{prob E-MIQP}, the statements of \cref{th MITR hardness,th MITR fesibility hardness} hold also for \cref{prob E-MIQP}.

We conclude the section by observing that theorem~1 in~\cite{dPDeyMol17MPA} implies that the problem of deciding whether the feasible region of \cref{prob E-MIQP} is nonempty lies in the complexity class NP.
Since \cref{prob MITR} can be mapped to \cref{prob E-MIQP} with the change of variables $y := B^{-1} (x - c)$, where $B$ is the invertible matrix with columns $b^1,\dots,b^n$, we obtain that also the problem of deciding whether the feasible region of \cref{prob MITR} is nonempty lies in the complexity class NP.
\cref{th MITR fesibility hardness} then implies that these two decision problems are NP-complete.


\section{Algorithmic and structural results for \cref{prob TR}}
\label{sec TR}

In this section we obtain a number of algorithmic and structural results for \cref{prob TR}
which we need to prove our main theorems.
In order to streamline the presentation, we introduce some notation.
A \emph{ball} with center $c$ and radius $r$ is a set of the form
\begin{align*}
\B(c,r) := \{x \in \R^n : \norm{x-c} \le r\},
\end{align*}
where $c \in \R^n$, and $r \in \R$ is positive. 
Notice that the feasible region of \cref{prob TR} 
is the ball $\B(0,1)$.
%
For ease of notation, from now on, we denote the objective function of \cref{prob TR} by 
$$
f(x) := x^\transp H x + h^\transp x.
$$
In this paper, we denote by $\log(x)$ the logarithm to the base $2$ of a positive $x \in \R$.

\subsection{Algorithmic results}

Vavasis and Zippel \cite{VavZip90} showed that we can find an arbitrarily good solution to \cref{prob TR} in polynomial time.
In this section, we discuss some algorithmic consequences of this result.




\begin{lemma}[\cite{VavZip90}]
\label{lem VavZip90}
Let $x^\circ$ be an optimal solution of \cref{prob TR}.
There is an algorithm which, for any positive integer $k$, finds a feasible solution $x^\dagger$ to \cref{prob TR} such that $$f(x^\dagger)-f(x^\circ) \le 2^{-k}.$$
The running time of the algorithm is polynomial in the size of \cref{prob TR} and in $k$.
\end{lemma}

\begin{prf}
Let $x^\dagger$ be the feasible solution to \cref{prob TR} obtained after $k'$ iterations with the algorithm in \cite{Ye92b,VavZip90}. 
From \cite{VavZip90} (see in particular the first three lines of section~5), we know $f(x^\dagger)-f(x^\circ) \le 2^{q-k'}$, where $q$ is a polynomial in $\size (H)$ and $\size (h)$, whose explicit form can be derived from section~4 in \cite{VavZip90}.
The result follows by setting $k':= \ceil{k+q}$, since each iteration of the algorithm performs a number of operations that is polynomial in $\size (H)$ and $\size (h)$.
\end{prf}

Next, we restate \cref{lem VavZip90} in such a way that the upper bound on the gap $f(x^\dagger) - f(x^\circ)$ is of the form $\epsilon \max \bra{\norm{H},\norm{h}}$.
Given a matrix $A \in \R^{m \times n}$, we denote by $\norm{A}$ its \emph{spectral norm}, which is defined by $\norm{A} := \sup_{x \neq 0} \norm{Ax} / \norm{x}$.
If $A$ is symmetric, $\norm{A}$ is equal to the maximum absolute value among the eigenvalues of $A$.

%
%
%
%
%
%

\begin{lemma}
\label{lem VavZip max}
Let $x^\circ$ be an optimal solution of \cref{prob TR}.
There is an algorithm which, for any positive $\epsilon \in \Q$, finds a feasible solution $x^\dagger$ to \cref{prob TR} such that $$f(x^\dagger)-f(x^\circ) \le \epsilon \max \bra{\norm{H},\norm{h}}.$$
The running time of the algorithm is polynomial in the size of \cref{prob TR} and in $\log(1/\epsilon)$.
\end{lemma}

\begin{prf}
If $H$ is the zero matrix and $h$ is the zero vector, we can simply set $x^\dagger := 0$.
Thus, in the remainder of the proof we assume that at least one among $H$ and $h$ is nonzero.
We set $k := \ceil{\log(1/\epsilon) + \size (H) + \size (h)}$, which can be computed in polynomial time.
From \cref{lem VavZip90}, it suffices to show
\begin{align*}
2^{-k} \le \epsilon \max \bra{\norm{H},\norm{h}}.
\end{align*}

We consider separately two cases.
In the first case, we assume $\norm{h} \ge \norm{H}$.
Then $h$ is nonzero, thus at least one component is nonzero, say $h_i$.
We obtain
$$
\norm{h}
= \sqrt{h_{1}^2 + \cdots + h_{n}^2}
\ge \abs{h_{i}} \ge 2^{-\size (h)}.
$$
Hence,
\begin{align*}
2^{-k} \le 2^{-\log(1/\epsilon) - \size (h)} = 2^{-\log(1/\epsilon)} 2^{-\size (h)} \le \epsilon \norm{h}
= \epsilon \max \bra{\norm{H},\norm{h}}.
\end{align*}

In the second case, we assume $\norm{H} \ge \norm{h}$.
Then $H$ is nonzero, thus at least one element is nonzero, say $h_{ij}$.
Let $e^j$ be the $j$th vector of the standard basis of $\R^n$.
We obtain
$$
\norm{H} = \sup_{x \neq 0} \frac{\norm{Hx}}{\norm{x}}
\ge \frac{\norm{H e^j}}{\norm{e^j}}
= \sqrt{h_{1j}^2 + \cdots + h_{nj}^2}
\ge \abs{h_{ij}} \ge 2^{-\size (H)}.
$$
Hence,
\begin{align*}
2^{-k} \le 2^{-\log(1/\epsilon) - \size (H)} = 2^{-\log(1/\epsilon)} 2^{-\size (H)} \le \epsilon \norm{H}
= \epsilon \max \bra{\norm{H},\norm{h}}.
\end{align*}
\end{prf}

%
%

We now show how \cref{lem VavZip max} can be used to compute a vector $v \in \B(0,1)$ with $|v^\transp H v|$ arbitrarily close to the spectral norm of $H$.


\begin{lemma}
\label{lem approx spectral norm multiplicative}
Let $H$ be a symmetric matrix in $\Q^{n \times n}$.
There is an algorithm which, for any positive $\epsilon \in \Q$, finds a vector $v \in \B(0,1)$ such that
$$
(1 - \epsilon) \norm{H} \le |v^\transp H v| \le \norm{H}.
$$ 
The running time of the algorithm is polynomial in $\size (H)$, $\log(1/\epsilon)$.
\end{lemma}

\begin{prf}
If $H$ is the zero matrix, we can simply set $v := 0$, thus in the remainder of the proof we assume $H$ nonzero.
Consider the following two optimization problems, which are special cases of \cref{prob TR}:
\begin{align*}
(S):= \ \min & \quad x^\transp H x & (T):= \ \max & \quad x^\transp H x & = \ - \min & \quad x^\transp (-H) x \\
\st & \quad x^\transp x \le 1, & \st & \quad x^\transp x \le 1 & \st & \quad x^\transp x \le 1.
\end{align*}
The algorithm is defined as follows.
We apply \cref{lem VavZip max} to both Problems $(S)$ and $(T)$, and 
find a feasible solution $s$ for $(S)$
and a feasible solution $t$ for $(T)$.
Return $v := t$ if $\abs{{t}^\transp H {t}} \ge \abs{{s}^\transp H {s}}$, and return $v := {s}$ otherwise.
In the remainder of the proof, we show the correctness of the algorithm.

Denote by $\lambda_{\min}$ the smallest eigenvalue of $H$ and by $\lambda_{\max}$ the largest eigenvalue of $H$.
Since $H$ is symmetric, we have $\lambda_{\min} = \min\bra{x^\transp H x : x^\transp x = 1}$ and $\lambda_{\max} = \max \bra{x^\transp H x : x^\transp x = 1}$.
Since $H$ is nonzero, we have $\norm{H} = \max \{|\lambda_{\min}|, \ |\lambda_{\max}|\} > 0$.
Every $x \in \R^n$ with $x^\transp x \le 1$ satisfies $|x^\transp H x| \le \norm{H}$, thus in the remainder of the proof it suffices to show $|v^\transp H v| \ge (1-\epsilon) \norm{H}$.
From the definition of our algorithm, it suffices to show either $|{s}^\transp H {s}| \ge (1-\epsilon) \norm{H}$ or $|{t}^\transp H {t}| \ge (1-\epsilon) \norm{H}$.
In the rest of the proof, we consider separately two cases.

In the first case, we assume $\norm{H} = |\lambda_{\min}|$.
We then have $\lambda_{\min} < 0$.
Hence, $\lambda_{\min} = \min \bra{ x^\transp H x : x^\transp x \le 1}$.
From \cref{lem VavZip max}, we know ${s}^\transp H s - \lambda_{\min} \le \epsilon \norm{H}$.
As $\lambda_{\min} = -\norm{H}$, this gives $|{s}^\transp H {s}| \ge -{s}^\transp H {s} \ge (1-\epsilon) \norm{H}$, as desired.

Next, we consider the case $\norm{H} = |\lambda_{\max}|$.
We then have $\lambda_{\max} > 0$.
Hence, $\lambda_{\max} = \max \bra{ x^\transp H x : x^\transp x \le 1}$.
From \cref{lem VavZip max}, we know $- {t}^\transp H t + \lambda_{\max} \le \epsilon \norm{H}$.
As $\lambda_{\max} = \norm{H}$, this gives  $|{t}^\transp H {t}| \ge {t}^\transp H {t} \ge (1-\epsilon) \norm{H}$, as desired.
\end{prf}



\subsection{Feasible points to \cref{prob TR} with distant objective values}

In this section, we provide two lower bounds on the difference, in objective value, between some specific feasible points to \cref{prob TR}.
First, we show that two vectors $s, t$, arbitrarily close to $h/\norm{h}$ and to $-h/\norm{h}$, respectively, yield a lower bound on $f(s)-f(t)$ that tends to $2 \norm{h}$, as $\epsilon$ goes to zero.
We remark that the bounds that we present are of interest when $\epsilon$ is close to zero.

\begin{lemma}
\label{lem linear candidate points}
Let $H$ be a symmetric matrix in $\R^{n \times n}$, let $h \in \R^n$,
and let $\epsilon \in [0,1/2]$.
Assume $\norm{H} \le \norm{h}$.
Let $l \in \R$ be positive such that $1-\epsilon \le l \norm{h} \le 1$.
Let $s \in \B((1-\epsilon)lh,\epsilon)$ and $t \in \B(-(1-\epsilon)lh,\epsilon)$.
Then
\begin{align*}
f(s)-f(t) \ge 2 \pare{1 - 5 \epsilon} \norm{h}.
\end{align*}
\end{lemma}

\begin{prf}
To show the result, we derive a lower bound for $f(s)$ and an upper bound for $f(t)$.
First, we give a lower bound for $f(s)$.
Let $g \in \R^n$ so that $s = (1-\epsilon)lh + g$.
Since $s \in \B((1-\epsilon)lh,\epsilon)$, we have $\norm{g} \le \epsilon$.
\begin{align*}
f(s) 
& =
{s}^\transp H s + h^\transp s \\
& =
((1-\epsilon)lh + g)^\transp H ((1-\epsilon)lh + g)  + h^\transp ((1-\epsilon)lh + g) \\
& =
(1-\epsilon)^2 l^2 h^\transp H h  + g^\transp H g + 2 (1-\epsilon) lh^\transp H g + (1-\epsilon) l \norm{h}^2 + h^\transp g \\
& \ge
(1-\epsilon)^2 l^2 h^\transp H h - \abs{g^\transp H g} - \abs{2 (1-\epsilon) lh^\transp H g} + (1-\epsilon) l \norm{h}^2 - \abs{h^\transp g} \\
& \ge
(1-\epsilon)^2 l^2 h^\transp H h - \norm{g}^2 \norm{H} - 2 (1-\epsilon) l \norm{h} \norm{H} \norm{g} + (1-\epsilon) l \norm{h}^2 - \norm{h} \norm{g} \\
& \ge
(1-\epsilon)^2 l^2 h^\transp H h - \epsilon^2 \norm{h} - 2 \epsilon (1-\epsilon) l \norm{h}^2 + (1-\epsilon) l \norm{h}^2 - \epsilon \norm{h} \\
& =
(1-\epsilon)^2 l^2 h^\transp H h - \epsilon^2 \norm{h} + (1- 2 \epsilon) (1-\epsilon) l \norm{h}^2 - \epsilon \norm{h} \\
& \ge
(1-\epsilon)^2 l^2 h^\transp H h - \epsilon^2 \norm{h} + (1- 2 \epsilon) (1-\epsilon)^2 \norm{h} - \epsilon \norm{h} \\
& =
(1-\epsilon)^2 l^2 h^\transp H h + \pare{1 - 5 \epsilon + 4 \epsilon^2 - 2 \epsilon^3} \norm{h}.
\end{align*}
Here, in the third inequality we used $\norm{g} \le \epsilon$, $\norm{H} \le \norm{h}$, and in the fourth inequality 
$l \norm{h} \ge 1-\epsilon$, since $1 - 2 \epsilon \ge 0$.

Similarly, we give an upper bound for $f(t)$.
Redefine $g \in \R^n$ so that $t = - (1-\epsilon)lh + g$.
Since $t \in \B(-(1-\epsilon)lh, \epsilon)$, we have $\norm{g} \le \epsilon$.
\begin{align*}
f(t) 
& =
{t}^\transp H t + h^\transp t \\
& =
(-(1-\epsilon)lh + g)^\transp H (-(1-\epsilon)lh + g)  + h^\transp (-(1-\epsilon)lh + g) \\
& =
(1-\epsilon)^2 l^2 h^\transp H h + g^\transp H g - 2 (1-\epsilon) lh^\transp H g - (1-\epsilon) l \norm{h}^2 + h^\transp g \\
& \le
(1-\epsilon)^2 l^2 h^\transp H h + \abs{g^\transp H g} + \abs{2 (1-\epsilon) lh^\transp H g} - (1-\epsilon) l \norm{h}^2 + \abs{h^\transp g} \\
& \le
(1-\epsilon)^2 l^2 h^\transp H h + \norm{g}^2 \norm{H} + 2 (1-\epsilon) l \norm{h} \norm{H} \norm{g} - (1-\epsilon) l \norm{h}^2 + \norm{h} \norm{g} \\
& \le
(1-\epsilon)^2 l^2 h^\transp H h + \epsilon^2 \norm{h} + 2 \epsilon (1-\epsilon) l \norm{h}^2 - (1-\epsilon) l \norm{h}^2 + \epsilon \norm{h} \\
& =
(1-\epsilon)^2 l^2 h^\transp H h + \epsilon^2 \norm{h} - (1-2 \epsilon) (1-\epsilon) l \norm{h}^2 + \epsilon \norm{h} \\
& \le
(1-\epsilon)^2 l^2 h^\transp H h + \epsilon^2 \norm{h} - (1-2 \epsilon) (1-\epsilon)^2 \norm{h} + \epsilon \norm{h} \\
& =
(1-\epsilon)^2 l^2 h^\transp H h - \pare{1 - 5 \epsilon + 4 \epsilon^2 - 2 \epsilon^3} \norm{h}.
\end{align*}
Here, in the third inequality we used 
$\norm{g} \le \epsilon$, $\norm{H} \le \norm{h}$, and in the fourth inequality $l \norm{h} \ge 1-\epsilon$, since $1 - 2 \epsilon \ge 0$.


We combine the two bounds above and obtain
\begin{align*}
f(s) - f(t)
& \ge
2 \pare{1 - 5 \epsilon + 4 \epsilon^2 - 2 \epsilon^3} \norm{h} \\
& \ge
2 \pare{1 - 5 \epsilon} \norm{h}.
\end{align*}
\end{prf}

Next, we show that two vectors $u, w$, arbitrarily close to the origin, and to a vector $v$ as in \cref{lem approx spectral norm multiplicative}, yield a lower bound on $\abs{f(w)-f(u)}$ that tends to $\norm{H}$, as $\epsilon$ goes to zero.

\begin{lemma}
\label{lem quadratic candidate points}
Let $H$ be a symmetric matrix in $\R^{n \times n}$, let $h \in \R^n$, and let $\epsilon \in [0,1]$.
Assume $\norm{h} \le \norm{H}$.
Let $v \in \B(0,1)$ such that $|v^\transp H v| \ge (1 - \epsilon) \norm{H}$.
Assume $h^\transp v \ge 0$ if $v^\transp H v \ge 0$, and $h^\transp v \le 0$ if $v^\transp H v<0$.
Let 
$w \in \B((1-\epsilon)v,\epsilon)$
and $u \in \B(0,\epsilon)$.
Then
$$
\abs{f(w)-f(u)} \ge \pare{1 - 7 \epsilon} \norm{H}.
$$
\end{lemma}

\begin{prf}
First, we consider the case $v^\transp H v \ge 0$.
To show the result, we derive an upper bound for $f(u)$ and a lower bound for $f(w)$.
We start with an upper bound for $f(u)$.
\begin{align*}
f(u) 
& =
{u}^\transp H u + h^\transp u \\
& \le 
\abs{{u}^\transp H u} + \abs{h^\transp u} \\
& \le 
\norm{u}^2 \norm{H} + \norm{h} \norm{u} \\
& \le 
\epsilon^{2} \norm{H} + \epsilon \norm{H} \\
& =
\pare{\epsilon + \epsilon^{2}} \norm{H}.
\end{align*}
Here, the third inequality holds because $\norm{u} \le \epsilon$ and $\norm{h} \le \norm{H}$.
%
Next, we give a lower bound for $f(w)$.
Let $g \in \R^n$ so that $w = (1-\epsilon)v + g$.
Since $w \in \B((1-\epsilon)v,\epsilon)$, we have $\norm{g} \le \epsilon$.
\begin{align*}
f(w) 
& =
{w}^\transp H w + h^\transp w \\
& = 
((1-\epsilon)v + g)^\transp H ((1-\epsilon)v + g) + h^\transp ((1-\epsilon)v + g) \\
& = 
(1-\epsilon)^2 v^\transp H v + 2 (1-\epsilon) v^\transp H g + g^\transp H g + (1-\epsilon) h^\transp v + h^\transp g \\
& \ge
(1-\epsilon)^2 v^\transp H v + 2 (1-\epsilon) v^\transp H g + g^\transp H g + h^\transp g \\
& \ge
(1-\epsilon)^2 v^\transp H v - \abs{2 (1-\epsilon) v^\transp H g} - \abs{g^\transp H g} - \abs{h^\transp g} \\
& \ge
(1-\epsilon)^3 \norm{H} - 2 (1-\epsilon) \norm{v} \norm{H} \norm{g} - \norm{g}^2 \norm{H} - \norm{h} \norm{g} \\
& \ge
(1-\epsilon)^3 \norm{H} - 2 \epsilon (1-\epsilon) \norm{H} - \epsilon^2 \norm{H} - \epsilon \norm{H} \\
& =
\pare{1 - 6 \epsilon + 4 \epsilon^2 - \epsilon^3} \norm{H}.
\end{align*}
Here, in the first inequality we used the assumption $h^\transp v \ge 0$, in the third inequality $v^\transp H v \ge (1 - \epsilon) \norm{H}$, and in the fourth inequality $\norm{v} \le 1$, $\norm{g} \le \epsilon$, $\norm{h} \le \norm{H}$.
%
We combine the two bounds above and obtain
\begin{align*}
f(w)-f(u)
& \geq
\pare{1 - 7 \epsilon + 3 \epsilon^2 - \epsilon^3} \norm{H} \\
& \geq
\pare{1 - 7 \epsilon} \norm{H}.
\end{align*}

\medskip

Next, we consider the case $v^\transp H v<0$.
This case is very similar to the previous one, and we provide a proof for completeness.
To show the result, we derive a lower bound for $f(u)$ and an upper bound for $f(w)$.
We start with a lower bound for $f(u)$.
\begin{align*}
f(u) 
& =
{u}^\transp H u + h^\transp u \\
& \ge 
- \abs{{u}^\transp H u} - \abs{h^\transp u} \\
& \ge 
- \norm{u}^2 \norm{H} - \norm{h} \norm{u} \\
& \ge 
- \epsilon^{2} \norm{H} - \epsilon \norm{H} \\
& =
- \pare{\epsilon + \epsilon^{2}} \norm{H}.
\end{align*}
Here, the third inequality holds because $\norm{u} \le \epsilon$ and $\norm{h} \le \norm{H}$.
%
Next, we give an upper bound for $f(w)$.
Let $g \in \R^n$ so that $w = (1-\epsilon)v + g$.
Since $w \in \B((1-\epsilon)v,\epsilon)$, we have $\norm{g} \le \epsilon$.
\begin{align*}
f(w) 
& =
{w}^\transp H w + h^\transp w \\
& = 
((1-\epsilon)v + g)^\transp H ((1-\epsilon)v + g) + h^\transp ((1-\epsilon)v + g) \\
& = 
(1-\epsilon)^2 v^\transp H v + 2 (1-\epsilon) v^\transp H g + g^\transp H g + (1-\epsilon) h^\transp v + h^\transp g \\
& \le
(1-\epsilon)^2 v^\transp H v + 2 (1-\epsilon) v^\transp H g + g^\transp H g + h^\transp g \\
& \le
(1-\epsilon)^2 v^\transp H v + \abs{2 (1-\epsilon) v^\transp H g} + \abs{g^\transp H g} + \abs{h^\transp g} \\
& \le
- (1-\epsilon)^3 \norm{H} + 2 (1-\epsilon) \norm{v} \norm{H} \norm{g} + \norm{g}^2 \norm{H} +  \norm{h} \norm{g} \\
& \le
- (1-\epsilon)^3 \norm{H} + 2 \epsilon (1-\epsilon) \norm{H} + \epsilon^2 \norm{H} + \epsilon \norm{H} \\
& =
- \pare{1 - 6 \epsilon + 4 \epsilon^2 - \epsilon^3} \norm{H}.
\end{align*}
Here, in the first inequality we used the assumption $h^\transp v \le 0$, in the third inequality $v^\transp H v \le - (1 - \epsilon) \norm{H}$, and in the fourth inequality $\norm{v} \le 1$, $\norm{g} \le \epsilon$, $\norm{h} \le \norm{H}$.
%
We combine the two bounds above and obtain
\begin{align*}
f(u)-f(w)
& \geq
\pare{1 - 7 \epsilon + 3 \epsilon^2 - \epsilon^3} \norm{H} \\
& \geq
\pare{1 - 7 \epsilon} \norm{H}.
\end{align*}
\end{prf}

\section{Feasible points to \cref{prob MITR} with distant objective values}
\label{sec MITR}

In this section we study \cref{prob MITR}.
The main goal is to introduce an algorithm that either finds two feasible solutions whose difference in objective value is at least $\max \bra{\norm{H},\norm{h}}$, for $\epsilon$ that goes to zero, or finds a nonzero vector $\direction$ along which the feasible region is flat.
This algorithm will be a key component in the design of the approximation algorithms for \cref{prob E-MIQP,prob MIQP} stated in \cref{th MITR algorithm,th MIQP algorithm}.
The reason we seek a bound featuring $\max \bra{\norm{H},\norm{h}}$ will become clear later, in the proof of \cref{prop E-MIQP flat}, and is essentially because the same $\max$ appears in the bound in \cref{lem VavZip max}.
%
To state the main result of this section, we need to introduce the concept of width.
Let $\S \subseteq \R^n$ be a nonempty bounded closed set and let $\direction \in \R^n$.
We define the \emph{width of~$\S$ along $\direction$} to be
\begin{align*}
\width_\direction (\S) := \max \bra{ \direction^\transp x : x \in \S } - \min \bra{ \direction^\transp x : x \in \S }.
\end{align*}

\begin{proposition}
[Distant points for \cref{prob MITR}]
\label{prop MITR dist points}
There is an algorithm which, for any $\epsilon \in (0,1/3]$,
either finds two feasible solutions $\badone,\badtwo$ to \cref{prob MITR} such that
\begin{align*}
f(\badone)-f(\badtwo) \geq 
\pare{1 - 7 \epsilon} \max \bra{\norm{H},\norm{h}},
\end{align*}
or finds a nonzero vector $\direction \in \Q^n$ with $\direction^\transp b^1, \dots, \direction^\transp b^p$ integer and $\direction^\transp b^{p+1} = \cdots = \direction^\transp b^n = 0$ such that $\width_\direction(\B(0,1)) \le  p \constlen/\epsilon$.
The running time of the algorithm is polynomial in the size of \cref{prob MITR} and in $\log(1/\epsilon)$.
\end{proposition}

\subsection{Ingredients}

In order to prove \cref{prop MITR dist points}, we need two more ingredients.
Our first ingredient concerns the approximation of square roots,
and is a consequence of \cref{lem sqrt approx}.
The key difference is that the approximation gap is multiplicative, rather than additive.

\begin{lemma}
\label{lem sqrt approx multiplicative}
Let $\alpha \in \Q$ be positive.
There is an algorithm which, for any positive $\epsilon \in \Q$, finds a positive $l \in \Q$ such that $l \le \sqrt{\alpha} \le (1+\epsilon) l$.
The running time of the algorithm is polynomial in $\size(\alpha), \log(1/\epsilon)$.
\end{lemma}

\begin{prf}
We use \cref{lem sqrt approx} with $k := \ceil{\log(1/\epsilon) + 3 \size (\alpha)/2+1}$.
We only need to show that the obtained $l,u$ satisfy $l > 0$ and $u \le (1+\epsilon) l$.

Using $k \ge 3 \size(\alpha)/2 + 1$ and $2^{-\size(\alpha)} \le \alpha \le 2^{\size(\alpha)}$, we obtain
\begin{align}
\label{lem sqrt approx multiplicative l bound}
\begin{split}
l
& = u - \frac{\max \bra{1,\alpha}}{2^k} 
\ge \sqrt{\alpha} - \frac{\max \bra{1,\alpha}}{2^{3 \size(\alpha)/2 + 1}} 
\ge 2^{-\size(\alpha)/2} - \frac{2^{\size(\alpha)}}{2^{3 \size(\alpha)/2 + 1}} \\
& 
= 2^{-\size(\alpha)/2} - 2^{-\size(\alpha)/2-1} 
= 2^{-\size(\alpha)/2-1},
\end{split}
\end{align}
which implies $l > 0$.

We now show $u \le (1+\epsilon) l$.
Using $k \ge \log(1/\epsilon) + 3 \size(\alpha)/2+1$, $\alpha \le 2^{\size (\alpha)}$, and \eqref{lem sqrt approx multiplicative l bound} we have
\begin{align*}
u & = l + \frac{\max \bra{1,\alpha}}{2^k} 
\le l + \frac{2^{\size(\alpha)}}{2^{\log(1/\epsilon) + 3 \size(\alpha)/2+1}} \\
& = l + 2^{-\size(\alpha)/2-1} 2^{\log(\epsilon)}
\le l + l \epsilon
= (1+\epsilon) l.
\end{align*}
\end{prf}

The second ingredient is a 
flatness result for balls, which follows with little effort from proposition~4 in \cite{dP23bMPA}.


\begin{lemma}[Flatness lemma]
\label{lem flat}
Let $a,c \in \Q^n$, let $\delta \in \Q$ with $\delta \ge 0$, let $b^1, \dots , b^n$ be linearly independent vectors in $\Q^n$, and let $p \in \{0,\dots,n\}$. 
There is a polynomial time algorithm which either finds a vector in $\B(a,\delta) \cap (\Pi_p(b^1, \dots , b^n) + \{c\})$, 
or finds a nonzero vector $\direction \in \Q^n$ with $\direction^\transp b^1, \dots, \direction^\transp b^p$ integer and $\direction^\transp b^{p+1} = \cdots = \direction^\transp b^n = 0$ such that $\width_\direction(\B(a,\delta)) \le p \constlen$.
\end{lemma}

\begin{prf}
In this proof, we denote by $\spn(\S)$ the linear span of a set $\S$ of vectors and by $\L^\perp$ the orthogonal complement of a linear subspace $\L$ or $\R^n$.
For $i=1,\dots,p$, let $c^i$ be the orthogonal projection of $b^i$ on $\spn(b^{p+1},\dots,b^n)^\perp$, i.e., $c^i := b^i - \sum_{j=p+1}^n ({b^j}^\transp b^i / \norm{b^j}^2) b^j$.
We apply proposition~4 in \cite{dP23bMPA} with $\Lambda := \Lambda(c^1,\dots,c^p)$, $d:=n$, and $a:=a-c$.
Then, there is a polynomial time algorithm which either finds a vector $\bar x \in \B(a-c,\delta) \cap (\Lambda+\spn(\Lambda)^\perp)$, or finds a nonzero vector $\direction \in \spn(\Lambda)$ with $\direction^\transp c^1, \dots, \direction^\transp c^p$ integer such that $\width_\direction(\B(a-c,\delta)) \le p \constlen$.

In the first case, the vector $\bar x' := \bar x + c$ is in $\B(a,\delta) \cap (\Lambda+\spn(\Lambda)^\perp + \{c\})$, and we have  $\Lambda+\spn(\Lambda)^\perp = \Lambda+ \spn(b^{p+1},\dots,b^n) = \Pi_p(b^1, \dots , b^n)$.
In the second case,
$\direction \in \spn(\Lambda)$ implies $\direction^\transp b^{p+1} = \cdots = \direction^\transp b^n = 0$.
We then have $\direction^\transp c^i = \direction^\transp b^i$, for $i=1,\dots,p$, thus 
$\direction^\transp c^1, \dots, \direction^\transp c^p$ integer implies $\direction^\transp b^1, \dots, \direction^\transp b^p$ integer.
Furthermore, by definition of width, we obtain $\width_\direction(\B(a,\delta)) = \width_\direction(\B(a-c,\delta)) \le p \constlen$.
\end{prf}

\subsection{Proof of \cref{prop MITR dist points}}



We are now ready to show \cref{prop MITR dist points}.
Before presenting its proof, we briefly discuss the 
the idea of the algorithm.
The goal is to either find two feasible solutions to \cref{prob MITR} whose difference in objective value is at least $\max \bra{\norm{H},\norm{h}}$, for $\epsilon$ that goes to zero, or to find a nonzero vector $\direction$ along which the feasible region is flat.
Depending on which one among $\norm{H}$ and $\norm{h}$ is largest, we search for the two feasible solutions in different regions in $\B(0,1)$;
As a result, we consider separately the cases $\norm{H} \ge \norm{h}$ and $\norm{h} > \norm{H}$.

In the case $\norm{H} \ge \norm{h}$, we would like to consider a vector $v \in \B(0,1)$ with $|v^\transp H v| = \norm{H}$.
However, this vector could be irrational, thus we use \cref{lem approx spectral norm multiplicative} to find, instead, a vector $v \in \B(0,1)$ with $|v^\transp H v|$ that is $\epsilon$-close to $\norm{H}$.
Next, we use our flatness lemma to seek two feasible points $u, w$ to \cref{prob MITR} that are $\epsilon$-close to the origin and to $v$, respectively.
If we find both these vectors, we use \cref{lem quadratic candidate points} to show that they have a difference in objective value that is at least $\norm{H}$, for $\epsilon$ that goes to zero.
Otherwise, we show that we have found a nonzero vector ($\direction$ in the statement) along which the feasible region $\B(0,1)$ is flat.

In the remaining case $\norm{h} > \norm{H}$, we would like to consider the vector $h/\norm{h}$.
However, $\norm{h}$ could be irrational, thus we find instead a scalar $l$ such that $l h \in \B(0,1)$ is $\epsilon$-close to $h/\norm{h}$.
Next, we use our flatness lemma to seek two feasible points $s, t$ to \cref{prob MITR} that are $\epsilon$-close to $lh$ and to $-lh$, respectively.
If we find both these vectors, we use \cref{lem linear candidate points} to show that they have a difference in objective value that is at least $\norm{h}$, for $\epsilon$ that goes to zero.
Otherwise, we show that we have found a nonzero vector ($\direction$ in the statement) along which the feasible region $\B(0,1)$ is flat.
We now give the complete proof of \cref{prop MITR dist points},
which uses \cref{lem sqrt approx multiplicative,lem approx spectral norm multiplicative,lem quadratic candidate points,lem linear candidate points,lem flat}.

\begin{prf}
In this proof, we consider separately the cases $\norm{H} \ge \norm{h}$ and $\norm{h} > \norm{H}$.

We start with the case $\norm{H} \ge \norm{h}$.
Let $\B_u := \B(0,\epsilon)$. 
From \cref{lem flat} (with $a:=0$ and $\delta:=\epsilon$), we either find a vector $u \in \B_u \cap (\Pi_p(b^1,\dots,b^n) + \{c\})$, or a nonzero vector $\direction \in \Q^n$ with $\direction^\transp b^1, \dots, \direction^\transp b^p$ integer and $\direction^\transp b^{p+1} = \cdots = \direction^\transp b^n = 0$
such that $\width_\direction(\B_u) \le p \constlen$. 
Using the algorithm in \cref{lem approx spectral norm multiplicative}, we find a vector $v \in \B(0,1)$ such that $|v^\transp H v| \ge (1 - \epsilon) \norm{H}$.
Note that the vector $-v$ also lies in $\B(0,1)$ and we have $(-v)^\transp H (-v) = v^\transp H v$.
Hence, by possibly redefining $v := -v$, we now assume $h^\transp v \ge 0$ if $v^\transp H v \ge 0$, and $h^\transp v \le 0$ if $v^\transp H v<0$.
Let $\B_w := \B((1-\epsilon)v,\epsilon)$.
From \cref{lem flat}, we either find a vector $w \in \B_w \cap (\Pi_p(b^1,\dots,b^n) + \{c\})$, or a nonzero vector $\direction \in \Q^n$ with $\direction^\transp b^1, \dots, \direction^\transp b^p$ integer and $\direction^\transp b^{p+1} = \cdots = \direction^\transp b^n = 0$ such that $\width_\direction(\B_w) \le p \constlen$.

Consider the subcase where \cref{lem flat} successfully found both vectors $u$ and $w$.
Note that 
$\B_u$ and $\B_w$ are contained in $\B(0,1)$,
thus $u, w$ are both feasible to \cref{prob MITR}.
We apply \cref{lem quadratic candidate points} and obtain
$$
\abs{f(w)-f(u)} \ge \pare{1 - 7 \epsilon} \norm{H}.
$$
%
%
Next, consider the subcase where \cref{lem flat} did not find both vectors $u$ and $w$.
In this case, \cref{lem flat} found a vector $\direction \in \Q^n$ with $\direction^\transp b^1, \dots, \direction^\transp b^p$ integer and $\direction^\transp b^{p+1} = \cdots = \direction^\transp b^n = 0$ such that $\width_\direction(\B) \le p \constlen$, for $\B \in \{\B_u, \B_w\}$. 
We then have
\begin{align*}
\width_\direction(\B(0,1)) 
& = \frac{1}{\epsilon} \width_\direction(\B) \le \frac{1}{\epsilon} p \constlen.
\end{align*}
This concludes the proof in the first case.

In the remainder of the proof, we consider the case $\norm{h} > \norm{H}$.
In particular, we have $h \neq 0$.
Using the algorithm in \cref{lem sqrt approx multiplicative} (with $\alpha:=1/\norm{h}^2$), we find a positive $l \in \Q$ such that $l \le 1/\norm{h} \le (1+\epsilon) l$.
We obtain $l \norm{h} \le 1$ and $l \norm{h} \ge 1 / (1+\epsilon) \ge 1-\epsilon$, thus $1-\epsilon \le l \norm{h} \le 1$.
Let $\B_s := \B((1-\epsilon)lh,\epsilon)$.
From \cref{lem flat}, we either find a vector $s \in \B_s \cap (\Pi_p(b^1,\dots,b^n) + \{c\})$, or a nonzero vector $\direction \in \Q^n$ with $\direction^\transp b^1, \dots, \direction^\transp b^p$ integer and $\direction^\transp b^{p+1} = \cdots = \direction^\transp b^n = 0$ such that $\width_\direction(\B_s) \le p \constlen$.
Let $\B_t := \B(-(1-\epsilon)lh,\epsilon)$.
From \cref{lem flat}, we either find a vector $t \in \B_t \cap (\Pi_p(b^1,\dots,b^n) + \{c\})$, or a nonzero vector $\direction \in \Q^n$ with $\direction^\transp b^1, \dots, \direction^\transp b^p$ integer and $\direction^\transp b^{p+1} = \cdots = \direction^\transp b^n = 0$ such that $\width_\direction(\B_t) \le p \constlen$.

Consider the subcase where \cref{lem flat} successfully found both vectors $s$ and $t$.
Note that
$\B_s$ and $\B_t$ are contained in $\B(0,1)$, thus $s, t$ are both feasible to \cref{prob MITR}.
%
We apply \cref{lem linear candidate points} and obtain 
$$
f(s)-f(t) \ge 2 \pare{1 - 5 \epsilon} \norm{h} \ge \pare{1 - 7 \epsilon} \norm{h},
$$
where the last inequality holds for $\epsilon \in (0,1/3]$.
%
%
Next, consider the subcase where \cref{lem flat} did not find both vectors $s$ and $t$.
In this case, \cref{lem flat} found a vector $\direction \in \Q^n$ with $\direction^\transp b^1, \dots, \direction^\transp b^p$ integer and $\direction^\transp b^{p+1} = \cdots = \direction^\transp b^n = 0$ such that $\width_\direction(\B) \le p \constlen$, for $\B \in \{\B_s, \B_t\}$. 
We then have
\begin{align*}
\width_\direction(\B(0,1)) 
& = \frac{1}{\epsilon} \width_\direction(\B) \le \frac{1}{\epsilon} p \constlen.
\end{align*}
This concludes the proof in the second case.
\end{prf}

\section{Flatness result for \cref{prob E-MIQP}}
\label{sec E-MIQP alt}

We are now ready to get back to \cref{prob E-MIQP}.
The main goal of this section is to present a flatness result for this problem:
an algorithm that either finds an $\epsilon$-approximate solution, or finds a nonzero vector $d$ along which the feasible region is flat.
%
To state this result, and to work with ellipsoids, it will be useful to introduce the following notation.
An \emph{ellipsoid} is a set of the form 
$$
\E(c,Q) := \bra{x \in \R^n : (x-c)^\transp Q (x-c) \leq 1},
$$
where $c \in \R^n$, and $Q \in \R^{n \times n}$ is a symmetric positive definite matrix. 
Note that $\B(c,r) = \E(c,I_n / r^2)$, where $I_n$ is the $n \times n$ identity matrix.
Notice that the inequality constraint in \cref{prob E-MIQP} 
can be equivalently written in the form $x \in \E(c,Q)$.
We say that an ellipsoid $\E(c,Q)$ is \emph{rational} when $c \in \Q^n$ and $Q \in \Q^{n \times n}$.
We can now state the main result of this section.

\begin{proposition}[Flatness result for \cref{prob E-MIQP}]
\label{prop E-MIQP flat}
There is an algorithm which, for any $\epsilon \in (0,1]$, either finds an $\epsilon$-approximate solution to \cref{prob E-MIQP}, 
or finds a nonzero vector $\direction \in \Z^n$ with $d_{p+1} = \cdots = d_n = 0$ and a scalar $\rho \in \Q$ such that 
\begin{align*}
\bra{d^\transp x : x \in \E(c,Q)} \subseteq \sbra{\rho, \rho + 19 p \constlen / \epsilon}.
\end{align*}
The running time of the algorithm is polynomial in the size of \cref{prob E-MIQP} and in $\log(1/\epsilon)$.
\end{proposition}

The key ingredient in the proof of \cref{prop E-MIQP flat} is \cref{prop MITR dist points}. 
However, we need
three additional results,
which are presented in \cref{sec E-MIQP alt ingredients}.

\subsection{Ingredients for flatness result}
\label{sec E-MIQP alt ingredients}

Given a rational ellipsoid $\E(c,Q)$, there may exist no rational affine transformation that maps $\E(c,Q)$ into the ball $\B(0,1)$.
An example is given by the ellipsoid $\E(0,(2))$, as we mentioned in \cref{sec intro}.
%
The goal of the next lemma is to show that for any positive $\delta \in \Q$, we can find in polynomial time a rational affine transformation that maps $\E(c,Q)$ into an ellipsoid that is arbitrarily close to the ball $\B(0,1)$.

\begin{lemma}
\label{lem sandwich ellipsoid balls}
Let $\E(c,Q) \subset \R^n$ be a rational ellipsoid.
There is an algorithm which, for any positive $\delta \in \Q$, finds a map $\tau: \R^n \to \R^n$ of the form $\tau(x) := B(x-c)$, with $B \in \Q^{n\times n}$ invertible, such that 
\begin{align*}
\B(0,1) \subseteq \tau(\E(c,Q)) \subseteq \B(0,1+\delta).
\end{align*}
The running time of the algorithm is polynomial in $\size(c), \size(Q), \log(1/\delta)$.
\end{lemma}

\begin{prf}
From corollary~1 in \cite{dP23bMPA}, 
there is a strongly polynomial algorithm that finds an invertible matrix
$M \in \Q^{n \times n}$ and a diagonal matrix $D \in \Q^{n \times n}$ such that $Q = M^\transp DM$.
Since $Q$ is positive definite, each diagonal element of $D$ is positive.
We can then apply \cref{lem sqrt approx multiplicative} to each diagonal element of $D$, and with $\epsilon := \delta$.
For $j = 1,\dots,n$, we consider the element $d_{jj}$ of $D$ in position $(j,j)$, and 
we obtain a positive $l_j \in \Q$ such that $l_j \le \sqrt{d_{jj}} \le (1+\delta) l_j$.
Let $L$ be the positive definite diagonal matrix in $\Q^{n \times n}$ with element $l_j$ in position $(j,j)$, for $j=1,\dots,n$.

Next, we show 
\begin{align}
\label{eq sandwich ellipsoid balls ineq y}
y^\transp y
\le 
y^\transp L^{-\transp} D L^{-1} y 
\le 
(1+\delta)^2 y^\transp y, \qquad \forall y \in \R^n.
\end{align}
To see this, we use $d_{jj} \ge l_j^2$ for $j=1,\dots,n$, to obtain
\begin{align*}
y^\transp L^{-\transp} D L^{-1} y 
= \sum_{j=1}^n \frac{d_{jj}}{l_j^2} y_j^2 
\ge \sum_{j=1}^n y_j^2 
= y^\transp y, \qquad \forall y \in \R^n,
\end{align*}
and we use $d_{jj} \le (1+\delta)^2 l_j^2$ for $j=1,\dots,n$, to derive
\begin{align*}
y^\transp L^{-\transp} D L^{-1} y 
= \sum_{j=1}^n \frac{d_{jj}}{l_j^2} y_j^2 
\le (1+\delta)^2 \sum_{j=1}^n y_j^2 
= (1+\delta)^2 y^\transp y, \qquad \forall y \in \R^n.
\end{align*}


Let $\tau: \R^n \to \R^n$ be the map defined by 
$\tau(x) := (1+\delta) LM(x-c)$,
thus with $B$ in the statement defined by $B := (1+\delta) LM$.
The containments in the statement then follow from \eqref{eq sandwich ellipsoid balls ineq y} by writing the corresponding three sets as follows:
\begin{align*}
\B(0,1) & = \bra{ y : y^\transp y \le 1}, \\
\tau(\E(c,Q)) & = \bra{\tau(x) : (x-c)^\transp Q (x-c) \le 1}, \\
& = \bra{y : y^{\transp} L^{-\transp} M^{-\transp} Q M^{-1} L^{-1} y \le (1+\delta)^2}, \\
& = \bra{y : y^\transp L^{-\transp} D L^{-1} y \le (1+\delta)^2}, \\
\B(0,1+\delta) & = \bra{y : y^\transp y \le (1+\delta)^2}.
\end{align*}
\end{prf}


The next result is a direct extension of \cref{lem VavZip max}, concerning the slightly different problem obtained from \cref{prob TR} by replacing the constraint $x \in \B(0,1)$ with $x \in \B(0,1+\delta)$:
\begin{align}
\label[problem]{prob TR delta}
\tag{TR$_\delta$}
\begin{split}
\min & \quad x^\transp H x + h^\transp x \\
\st & \quad x^\transp x \le (1+\delta)^2.
\end{split}
\end{align}
Here $H$ is a symmetric matrix in $\Q^{n \times n}$, $h \in \Q^n$, and $\delta \in \Q$ is nonnegative.


\begin{lemma}
\label{lem VavZip max delta}
Let $x^\circ$ be an optimal solution of \cref{prob TR delta}.
There is an algorithm which, for any positive $\epsilon \in \Q$, finds a feasible solution $x^\dagger$ to \cref{prob TR delta} such that $$f(x^\dagger)-f(x^\circ) \le \epsilon (1+\delta)^2 \max \bra{\norm{H},\norm{h}}.$$
The running time of the algorithm is polynomial in the size of \cref{prob TR delta} and in $\log(1/\epsilon)$.
\end{lemma}

\begin{prf}
We apply the change of variables $y := x / (1+\delta)$ to \cref{prob TR delta} and obtain
\begin{align}
\label[problem]{prob TR delta prime}
\tag{TR$'_\delta$}
\begin{split}
\min & \quad y^\transp (1+\delta)^2 H y + (1+\delta) h^\transp y \\
\st & \quad y^\transp y \le 1.
\end{split}
\end{align}
For ease of notation, we denote by $f_\delta(y)$ the objective function of \cref{prob TR delta prime}. 
Note that $y^\circ := x^\circ / (1+\delta)$ is an optimal solution of \cref{prob TR delta prime}.
Let $y^\dagger$ be a feasible solution to \cref{prob TR delta prime} obtained with \cref{lem VavZip max}.
We let $x^\dagger := (1+\delta) y^\dagger$ and note that $x^\dagger$ is feasible to \cref{prob TR delta}.
We obtain 
\begin{align*}
f(x^\dagger)-f(x^\circ) & = f_\delta(y^\dagger)-f_\delta(y^\circ) \\
& \le \epsilon \max \bra{(1+\delta)^2\norm{H},(1+\delta)\norm{h}} \\
& \le \epsilon (1+\delta)^2 \max \bra{\norm{H},\norm{h}}.
\end{align*}
\end{prf}

The last result of the section discusses how flat directions in a mixed integer lattice are mapped to flat directions in the mixed integer lattice $\Z^p \times \R^{n-p}$.



\begin{lemma}
\label{lem direction transformation}
Let $b^1,\dots,b^n$ be linearly independent vectors in $\R^n$ and let $p \in \{0,\dots,n\}$.
Let $B \in \R^{n\times n}$ with columns $b^1,\dots,b^n$, let $c \in \R^n$, and let $\tau : \R^n \to \R^n$ be the map defined by $\tau(x) := B(x-c)$.
Let ${\direction'} \in \R^n$ with ${\direction'}^\transp b^1, \dots, {\direction'}^\transp b^p$ integer and ${\direction'}^\transp b^{p+1} = \cdots = {\direction'}^\transp b^n = 0$.
Define $\direction := B^\transp {\direction'}$.
Then, $d \in \Z^n$, $d_{p+1} = \cdots = d_n = 0$, and, for any nonempty bounded closed set $\S \subseteq \R^n$, $\width_\direction(\S) = \width_{\direction'}(\tau(\S)).$
\end{lemma}

\begin{prf}
For $j=1,\dots,n$, we have $d_j = {b^j}^\transp d'$.
Thus, ${\direction'}^\transp b^1, \dots, {\direction'}^\transp b^p$ integer implies $d_1, \cdots, d_p$ integer, and ${\direction'}^\transp b^{p+1} = \cdots = {\direction'}^\transp b^n = 0$ implies $d_{p+1} = \cdots = d_n = 0$.
We now show that for any nonempty bounded closed set $\S \subseteq \R^n$, we have 
$
\width_\direction(\S) = \width_{\direction'}(\tau(\S)).
$
\begin{align*}
\width_\direction(\S)
& = \max \bra{ \direction^\transp x : x \in \S } - \min \bra{ \direction^\transp x : x \in \S } \\
& = \max \bra{ \direction^\transp (x-c) : x \in \S } - \min \bra{ \direction^\transp (x-c) : x \in \S } \\
& = \max \bra{ {\direction'}^\transp B (x-c) : x \in \S } - \min \bra{ {\direction'}^\transp B (x-c) : x \in \S } \\
& = \max \bra{ {\direction'}^\transp y : y \in \tau(\S) } - \min \bra{ {\direction'}^\transp y : y \in \tau(\S) } \\
& = \width_{\direction'}(\tau(\S)).
\end{align*}
%
\end{prf}

\subsection{Proof of \cref{prop E-MIQP flat}}
\label{sec MITR flat}


We now have all ingredients that we need to prove \cref{prop E-MIQP flat}.
Before giving the proof, we give an overview of the algorithm.
First, we construct an affine transformation $\tau$ that maps $\E(c,Q)$ into an ellipsoid that is sandwiched between $\B(0,1)$ and $\B(0,1+\epsilon)$.
We denote by \cref{prob E-MIQP'} the optimization problem obtained by applying the same affine transformation to \cref{prob E-MIQP}.
Next, we consider \cref{prob E-MIQP'-W} obtained from \cref{prob E-MIQP'} by replacing the constraint $y \in \tau( \E(c,Q) )$ with the weaker requirement $y \in \B(0,1 + \epsilon)$ and by dropping the mixed integer lattice constraint.
We then find an almost optimal solution $y^\dagger$ to \cref{prob E-MIQP'-W} and use \cref{lem flat} to search for a close by vector $y^\diamond$ that is feasible to \cref{prob E-MIQP'}.
Next, we consider \cref{prob E-MIQP'-S} obtained from \cref{prob E-MIQP'} by replacing the constraint $y \in \tau( \E(c,Q) )$ with with the stronger requirement $y \in \B(0,1)$.
Using \cref{prop MITR dist points}, we search for two feasible solutions $\badoney,\badtwoy$ to \cref{prob E-MIQP'-S} that are distant in objective value.
If we manage to find all three vectors $y^\diamond$, $\badoney$, $\badtwoy$, we show that $y^\diamond$ is an $\epsilon$-approximate solution to \cref{prob E-MIQP'}.
Otherwise, we show that $\E(c,Q)$ is flat and we find the associated vector $d$.
%
We now give the complete proof of \cref{prop E-MIQP flat}, which uses \cref{lem sandwich ellipsoid balls,lem VavZip max delta,lem flat,prop MITR dist points,lem direction transformation}.

\begin{prf}
Let $\epsilon' := \epsilon/18 \in (0,1/18]$. 
Using the algorithm in \cref{lem sandwich ellipsoid balls} (with $\delta := \epsilon'$), we find a map $\tau: \R^n \to \R^n$ of the form $\tau(x) := B(x-c)$, with $B \in \Q^{n\times n}$ invertible, such that 
\begin{align*}
\B(0,1) \subseteq \tau( \E(c,Q) ) \subseteq \B(0,1+\epsilon').
\end{align*} 
%
%
Applying the change of variables $y := B (x-c)$ to \cref{prob E-MIQP}, we obtain the optimization problem
\begin{align}
\label[problem]{prob E-MIQP'}
\tag{E-MIQP$'$}
\begin{split}
\gamma
+ 
\min & \quad y^\transp H' y 
+ {h'}^\transp y
 \\
\st & \quad y \in \tau( \E(c,Q) ) \\
& \quad y \in 
\Pi_p(b^1,\dots,b^n) + \{c'\},
\end{split}
\end{align}
where $\gamma:= c^\transp H c + {h}^\transp c$, $H' := B^{-\transp} H B^{-1}$, 
$h' := 2 B^{-\transp} H^\transp c + B^{-\transp} h$, 
$b^1,\dots,b^n$ are the columns of $B$,
and $c':=-Bc$.
For ease of notation, in the following we denote by $f_y(y) := y^\transp H' y + {h'}^\transp y$.


Consider the optimization problem obtained from \cref{prob E-MIQP'} by replacing the constraint $y \in \tau( \E(c,Q) )$ with the weaker requirement $y \in \B(0,1 + \epsilon')$ and by dropping the mixed integer lattice constraint:
\begin{align}
\label[problem]{prob E-MIQP'-W}
\tag{E-MIQP$'$-w}
\begin{split}
\gamma + 
\min & \quad y^\transp H' y 
+ {h'}^\transp y
 \\
\st & 
\quad y^\transp y \le (1 + \epsilon')^2.
\end{split}
\end{align}
Let $y^\circ$ be an optimal solution of \cref{prob E-MIQP'-W}, which we do not need to compute.
Using the algorithm in \cref{lem VavZip max delta} (with $\epsilon := \delta := \epsilon'$), we find a feasible solution $y^\dagger$ to \cref{prob E-MIQP'-W} such that 
\begin{align}
\label{eq MITR from lem VavZip new delta}
f_y(y^\dagger)-f_y(y^\circ) \le \epsilon' (1+\epsilon')^2 \max \bra{\norm{H'},\norm{h'}}.
\end{align}
Let $\B_\diamond := \B((1-\epsilon')y^\dagger / (1+\epsilon'),\epsilon')$. 
From \cref{lem flat}, we either find a vector $y^\diamond \in \B_\diamond \cap (\Pi_p(b^1,\dots,b^n) + \{c'\})$, 
or a nonzero vector ${\direction'} \in \Q^n$ with ${\direction'}^\transp b^1, \dots, {\direction'}^\transp b^p$ integer and ${\direction'}^\transp b^{p+1} = \cdots = {\direction'}^\transp b^n = 0$ such that
$\width_{\direction'}(\B_\diamond) \le p \constlen$. 

Consider now the mixed integer trust region problem obtained from \cref{prob E-MIQP'} by replacing the constraint $y \in \tau( \E(c,Q) )$ with the stronger requirement $y \in \B(0,1)$:
\begin{align}
\label[problem]{prob E-MIQP'-S}
\tag{E-MIQP$'$-s}
\begin{split}
\gamma + 
\min & \quad y^\transp H' y 
+ {h'}^\transp y
 \\
\st & \quad y^\transp y \le 1 \\
& \quad y \in \Pi_p(b^1,\dots,b^n) + \{c'\}.
\end{split}
\end{align}
Using the algorithm in \cref{prop MITR dist points}, we either find two feasible solutions $\badoney,\badtwoy$ to \cref{prob E-MIQP'-S} such that
\begin{align}
\label{eq E-MIQP from prop MITR dist points}
f(\badoney)-f(\badtwoy) \geq 
\pare{1 - 7 \epsilon'} \max \bra{\norm{H'},\norm{h'}},
\end{align}
or we find a nonzero vector
${\direction'} \in \Q^n$ with ${\direction'}^\transp b^1, \dots, {\direction'}^\transp b^p$ integer and ${\direction'}^\transp b^{p+1} = \cdots = {\direction'}^\transp b^n = 0$ such that $\width_{\direction'}(\B(0,1)) \le  p \constlen/\epsilon'$.

In the remainder of the proof, we consider separately two cases.
First, consider the case where \cref{lem flat,prop MITR dist points} successfully found all three vectors $y^\diamond$, $\badoney$, $\badtwoy$.
Note that, since $y^\dagger \in \B(0,1+\epsilon')$, we have $y^\dagger/(1+\epsilon') \in \B(0,1)$, thus $\B_\diamond \subseteq \B(0,1)$ and
$y^\diamond$ is feasible to \cref{prob E-MIQP'-S}.
Hence, the three vectors $y^\diamond$, $\badoney$, $\badtwoy$ are feasible to \cref{prob E-MIQP'-S} and to \cref{prob E-MIQP'}.
%
%
Denote by $f_{\sup}$ and $f_{\inf}$ the supremum and the infimum of $f_y(y)$ on the feasible region of \cref{prob E-MIQP'}.
From \eqref{eq E-MIQP from prop MITR dist points}, we directly obtain the following lower bound on $f_{\sup}-f_{\inf}$:
\begin{align}
\label{eq MITR lower bound}
\begin{split}
f_{\sup}-f_{\inf} & \ge f(\badoney)-f(\badtwoy) \\
& \ge \pare{1 - 7 \epsilon'} \max \bra{\norm{H'},\norm{h'}}.
\end{split}
\end{align}
%
Next, we obtain an upper bound on $f_y(y^\diamond)-f_{\inf}$.
First, let $g \in \R^n$ so that $y^\diamond = y^\dagger + g$. 
It is simple to check that $\B_\diamond \subseteq \B(y^\dagger,3\epsilon')$, which implies $\norm{g} \le 3\epsilon'$.
Then,
\begin{align*}
\abs{f_y(y^\diamond)-f_y(y^\dagger)}
& = 
\abs{f_y(y^\diamond)-f_y(y^\diamond-g)} \\
& =
\abs{{y^\diamond}^\transp H' y^\diamond + {h'}^\transp y^\diamond - (y^\diamond - g)^\transp H' (y^\diamond - g) - {h'}^\transp (y^\diamond - g)} \\
& = 
\abs{2 {y^\diamond}^\transp H' g - g^\transp H' g + {h'}^\transp g} \\
& \le
\abs{2 {y^\diamond}^\transp H' g} + \abs{g^\transp H' g} + \abs{{h'}^\transp g} \\
& \le
2 \norm{y^\diamond} \norm{H'} \norm{g} + \norm{g}^2 \norm{H'} + \norm{{h'}} \norm{g} \\
& \le
6 \epsilon' \norm{H'} + 9 {\epsilon'}^2 \norm{H'} + 3 \epsilon' \norm{{h'}} \\
& \le
\pare{9 \epsilon' + 9 {\epsilon'}^2} \max \bra{\norm{H'}, \norm{h'}}.
\end{align*}
Here, in the third inequality we used $\norm{y^\diamond} \le 1$, $\norm{g} \le 3\epsilon'$.
Using the above relation and \eqref{eq MITR from lem VavZip new delta}, we obtain the following upper bound on $f_y(y^\diamond)-f_{\inf}$:
\begin{align}
\label{eq MITR upper bound}
\begin{split}
f_y(y^\diamond)-f_{\inf} 
& \le f_y(y^\diamond)-f_y(y^\circ) \\
& \le \abs{f_y(y^\diamond)-f_y(y^\dagger)} + (f_y(y^\dagger)-f_y(y^\circ)) \\
& \le \pare{9 \epsilon' + 9 {\epsilon'}^2 + \epsilon' (1+\epsilon')^2} \max \bra{\norm{H'},\norm{h'}} \\
& = \pare{10 \epsilon' + 11 \epsilon'^2 + \epsilon'^3} \max \bra{\norm{H'},\norm{h'}}.
\end{split}
\end{align}

We are now ready to show that $y^\diamond$ is an $\epsilon$-approximate solution to \cref{prob E-MIQP'}.
We have
\begin{align*}
f(y^\diamond)-f_{\inf}
& \le (10 \epsilon' + 11 \epsilon'^2 + \epsilon'^3) \max \bra{\norm{H'},\norm{h'}} \\
& \le 18 \epsilon' (1 - 7 \epsilon') \max \bra{\norm{H'},\norm{h'}} \\
& \le 18 \epsilon' (f_{\sup}-f_{\inf}) \\
& = \epsilon (f_{\sup}-f_{\inf}).
\end{align*}
Here, in the first inequality we used \eqref{eq MITR upper bound},
the second inequality holds for $\epsilon' \in (0, 1/18]$,
and in the third inequality we used \eqref{eq MITR lower bound}.
It then follows that $x^\diamond := B^{-1} y^\diamond + c$ is an $\epsilon$-approximate solution to \cref{prob E-MIQP}.
This concludes the first case.

Next, consider the case where 
\cref{lem flat} did not find the vector $y^\diamond$, or 
\cref{prop MITR dist points} did not find the two vectors $\badoney$, $\badtwoy$.
In this case, we first show that we have found a nonzero vector ${\direction'} \in \Q^n$ with ${\direction'}^\transp b^1, \dots, {\direction'}^\transp b^p$ integer and ${\direction'}^\transp b^{p+1} = \cdots = {\direction'}^\transp b^n = 0$ such that
\begin{align*}
\width_{\direction'}(\B(0,1)) \le  \frac{1}{\epsilon'} p \constlen.
\end{align*}
If \cref{prop MITR dist points} did not find the two vectors $\badoney$, $\badtwoy$, then we already know it.
Otherwise, if \cref{lem flat} did not find the vector $y^\diamond$, then we have found a nonzero vector ${\direction'} \in \Q^n$ with ${\direction'}^\transp b^1, \dots, {\direction'}^\transp b^p$ integer and ${\direction'}^\transp b^{p+1} = \cdots = {\direction'}^\transp b^n = 0$ such that $\width_{\direction'}(\B_\diamond) \le p \constlen$. 
We obtain
\begin{align*}
\width_{\direction'}\pare{\B(0,1)}
& = \frac{1}{\epsilon'} \width_{\direction'}\pare{\B_\diamond} 
\le \frac{1}{\epsilon'} p \constlen.
\end{align*}

We can now upper bound $\width_{\direction'}\pare{\B(0,1+\epsilon')}$ as follows:
\begin{align*}
\width_{\direction'}\pare{\B(0,1+\epsilon')} 
& = (1+\epsilon') \width_{\direction'}\pare{\B(0,1)} \\
& \le \frac{(1+\epsilon')}{\epsilon'} p \constlen \\
& = \pare{1+\frac{18}{\epsilon}} p \constlen \\
& \le \frac{19}{\epsilon} p \constlen.
\end{align*}


Now let $\direction := B^\transp {\direction'} \in \Q^n$.
From \cref{lem direction transformation}, $d \in \Z^n$, $d_{p+1} = \cdots = d_n = 0$, and, for any nonempty bounded closed set $\S \subseteq \R^n$, $\width_\direction(\S) = \width_{\direction'}(\tau(\S)).$
If we let $\S := \tau^\leftarrow(\B(0,1+\epsilon'))$, where $\tau^\leftarrow$ denotes the inverse of $\tau$, we obtain
\begin{align*}
\width_\direction(\tau^\leftarrow(\B(0,1+\epsilon'))) 
= \width_{\direction'}(\B(0,1+\epsilon'))
\le \frac{19}{\epsilon} p \constlen.
\end{align*}
%
Let $\rho := \min\{{\direction}^\transp x : x \in \tau^\leftarrow(\B(0,1+\epsilon'))\}$, and note that we can calculate it as follows:
\begin{align*}
\rho 
& = \min\{{\direction}^\transp (B^{-1} y + c) : y \in \B(0,1+\epsilon')\} \\
& = \min\{{\direction'}^\transp y : y \in \B(0,1+\epsilon')\} + \direction^\transp c \\
& = - (1+\epsilon') \norm{\direction'} + \direction^\transp c \\
& = - (1+\epsilon/18) \norm{\direction'} + \direction^\transp c.
\end{align*}
We obtain
\begin{align*}
\{d^\transp x : x \in \tau^\leftarrow(\B(0,1+\epsilon'))\} \subseteq [\rho, \rho + 19 p \constlen / \epsilon].
\end{align*}
Since $\E(c,Q) \subseteq \tau^\leftarrow(\B(0,1+\epsilon'))$, we also have 
\begin{align*}
\{d^\transp x : x \in \E(c,Q)\} \subseteq [\rho, \rho + 19 p \constlen / \epsilon].
\end{align*}
%
This concludes the second case.
\end{prf}

\section{Approximation algorithm for \cref{prob E-MIQP}}
\label{sec E-MIQP alg}

In this section we prove \cref{th MITR algorithm}.
Before diving into the proof, we present, in \cref{sec E-MIQP alg ingredients}, a few additional results that we need.

\subsection{Ingredients for approximation algorithm}
\label{sec E-MIQP alg ingredients}


We start by presenting two lemmas that will be instrumental in reformulating a restricted \cref{prob E-MIQP}, where the feasible region is intersected with a hyperplane, as a new \cref{prob E-MIQP} with one fewer integer variable.
In the first lemma we show that, given a hyperplane $\H$ in $\R^n$ that satisfies some technical conditions, we can find an affine transformation that maps the hyperplane $\{y \in \R^n : y_1 = 0\}$ to $\H$ and that preserves mixed integer points.
\cref{lem affine transformation} is similar to lemma~6 in \cite{dP23bMPA}.
The main difference is that here we seek an affine function from $\R^n$ to $\R^n$, while in lemma~6 in \cite{dP23bMPA} the affine function is from $\R^{n-1}$ to $\H$.
It will be useful to introduce some notation for hyperplanes.
A \emph{hyperplane} is a set of the form 
$$
\H(\direction,\beta) := \{x \in \R^n : \direction^\transp x = \beta \},
$$
where $\direction \in \R^n \setminus \{0\}$, $\beta \in \R$.
We say that a hyperplane is \emph{rational}, when $\direction \in \Q^n$ and $\beta \in \Q$.
Recall that a square matrix $T \in \Z^{n\times n}$ is called \emph{unimodular} if it has determinant $\pm 1$.
The arguments used in the proof are direct extensions of those for pure integer linear programs (see, e.g., \cite{ConCorZamBook}).
%



\begin{lemma}
\label{lem affine transformation}
Let $\H(\direction,\beta) \subset \R^n$ be a rational hyperplane with $d_{p+1} = \cdots = d_n = 0$, for $p \in \{1,\dots,n\}$.
There is a polynomial time algorithm which determines whether the set $\S:= \H(\direction,\beta) \cap (\Z^p \times \R^{n-p})$ is empty or not.
If $\S \neq \emptyset$, the algorithm finds a map $\eta: \R^n \to \R^n$ of the form $\eta(y) := \bar x + Ty$, with $\bar x \in \Z^n$ and $T \in \Z^{n\times n}$ unimodular, such that
\begin{align*}
\H(\direction,\beta) & = \eta ( \{y \in \R^n : y_1=0 \} ),
\\
\Z^p \times \R^{n-p} & = \eta( \Z^p \times \R^{n-p} ).
\end{align*}
\end{lemma}

\begin{prf}
%
By possibly multiplying the equation $\direction^\transp x = \beta$ by the least common multiple of the denominators of the entries of $\direction$, we may assume that $\direction$ is an integral vector.
By possibly dividing the equation $\direction^\transp x = \beta$ by the greatest common divisor of the entries of $\direction$, we may assume that $\direction$ has relatively prime entries.
If $\beta \notin \Z$, then $\S$ is empty and we are done.
Thus, we now assume $\beta \in \Z$.
Since $\direction_1,\dots,\direction_p$ are relatively prime, by corollary~1.9 in~\cite{ConCorZamBook}, the equation $\sum_{j=1}^p \direction_j x_j = \beta$ has an integral solution, thus $\S$ is nonempty.
It is well known (see, e.g., proof of corollary 1.9 in~\cite{ConCorZamBook}) that we can compute in polynomial time a vector $\tilde x \in \Z^p$ and a unimodular matrix $U \in \Z^{p \times p}$ such that 
\begin{align*}
\left\{x \in \Z^p : \sum_{j=1}^p \direction_j x_j= \beta\right\} = \{\tilde x + Uy : y \in \Z^{p}, \ y_1=0\}.
\end{align*}
Taking the convex hull of the two sets in the above equality, we obtain
\begin{align*}
\left\{x \in \R^p : \sum_{j=1}^p \direction_j x_j= \beta\right\} = \{\tilde x + Uy : y \in \R^{p}, \ y_1=0\}.
\end{align*}
Since $\tilde x \in \Z^p$ and $U$ is unimodular we also have
\begin{align*}
\Z^p = \{\tilde x + Uy : y \in \Z^{p}\}.
\end{align*}

We define the vector $\bar x \in \Z^n$ by $\bar x_j := \tilde x_j$ for $j = 1,\dots,p$, and $\bar x_j := 0$ for $j = p+1,\dots,n$. 
We also define the unimodular matrix $T \in \Z^{n \times n}$ with block corresponding to the first $p$ rows and $p$ columns being equal to $U$, block corresponding to the last $n-p$ rows and $n-p$ columns being equal to the identity matrix $I_{n-p}$, and remaining entries zero, i.e.,
$$
T:=\left(\begin{array}{c|c}
U & 0 \\
\hline 
0 & I_{n-p}
\end{array}\right). 
$$
Since $\direction_i = 0$ for all $i \in \{p+1, \dots,n\}$, we conclude 
\begin{align*}
\H(\direction,\beta) & = \{\bar x + Ty : y \in \R^n, \ y_1=0 \} \\
\Z^p \times \R^{n-p} & = \{\bar x + Ty : y \in \Z^p \times \R^{n-p}\}.
\end{align*}
The result then follows by defining the rational affine transformation $\eta: \R^n \to \R^n$, 
as $\eta(y) := \bar x + Ty$.
\end{prf}

The next lemma builds on \cref{lem affine transformation} and is of key importance in our path to prove \cref{th MITR algorithm}.
Given a hyperplane $\H$ and an ellipsoid $\E$ in $\R^n$, we show that we can find an ellipsoid $\E'$ in $\R^{n-1}$ and an affine transformation from $\R^{n-1}$ to $\R^n$, which maps $\E'$ to $\E \cap \H$ and that preserves mixed integer points.

\begin{lemma}
\label{lem mitr go to lower dim}
Let $\E(c,Q) \subset \R^n$ be a rational ellipsoid and let $\H(\direction,\beta) \subset \R^n$ be a rational hyperplane with $d_{p+1} = \cdots = d_n = 0$, for $p \in \{1,\dots,n\}$.
There is a polynomial time algorithm which determines whether the sets $\S:=\H(\direction,\beta) \cap ( \Z^p \times \R^{n-p} )$ and $\E(c,Q) \cap \H(\direction,\beta)$ are empty or not.
If $\S \neq \emptyset$ and $n \ge 2$, the algorithm finds 
a map $\eta: \R^{n-1} \to \R^n$ of the form $\eta(y) := \bar x + Ty$, with $\bar x \in \Z^n$ and $T \in \Q^{n\times {n-1}}$ of full rank, such that
\begin{align*}
\H(\direction,\beta) & =  \eta(\R^{n-1}), \\
\S & =
\eta(\Z^{p-1} \times \R^{n-p}).
\end{align*} 
If also $\E(c,Q) \cap \H(\direction,\beta) \neq \emptyset$, the algorithm finds the preimage of the set $\E(c,Q) \cap \H(\direction,\beta)$ under $\eta$, which is either a rational ellipsoid $\E(c',Q') \subset \R^{n-1}$ or a singleton $\bra{c'}$, for some $c' \in \Q^{n-1}$.
\end{lemma}

\begin{prf}
%
From \cref{lem affine transformation},
there is a polynomial time algorithm that determines whether the set $\S$ is empty or not.
If $\S \neq \emptyset$, the algorithm finds a map $\eta': \R^n \to \R^n$ of the form $\eta'(y) := \bar x + T'y$, with $\bar x \in \Z^n$ and $T' \in \Z^{n\times n}$ unimodular, such that
\begin{align*}
\H(\direction,\beta) & = \eta' ( \{y \in \R^n : y_1=0 \} ) \\
\Z^p \times \R^{n-p} & = \eta'( \Z^p \times \R^{n-p} ).
\end{align*}
%
%
The preimage of the ellipsoid $\E(c, Q)$ under $\eta'$ is the ellipsoid 
\begin{align*}
& \bra{y \in \R^n : (\bar x + T'y - c)^\transp Q (\bar x + T'y - c) \leq 1} = \\
=
&\bra{y \in \R^n : (y - ({T'}^{-1} (c- \bar x)))^\transp {T'}^\transp Q {T'} (y - ({T'}^{-1} (c- \bar x))) \leq 1} \\
= 
& \ \E({T'}^{-1} (c- \bar x),{T'}^\transp Q {T'}).
\end{align*}
For ease of notation, we define $q:={T'}^{-1} (c- \bar x)$ and $M := {T'}^\transp Q {T'}$, 
so that the preimage of $\E(c, Q)$ under $\eta'$ is $\E(q,M)$.


Let $T \in \Q^{n \times n-1}$ be the matrix $T'$ stripped of the first column, and define the map $\eta: \R^{n-1} \to \R^n$ by $\eta(y) := \bar x + Ty$.
We obtain $\H(\direction,\beta) =  \eta(\R^{n-1})$ and
$
\S =
\eta(\Z^{p-1} \times \R^{n-p}).
$
%
Let $\E' \subset \R^{n-1}$ be the preimage of the set $\E(c,Q) \cap \H(\direction,\beta)$ under $\eta$.
We then have
\begin{align*}
\E' := \bra{\tilde y \in \R^{n-1} : (0,\tilde y) \in \E(q,M) }.
\end{align*}
To conclude the proof, we need to show that $\E'$ is either a rational ellipsoid $\E(c',Q') \subset \R^{n-1}$ or a singleton $\bra{c'}$, for some $c' \in \Q^{n-1}$.
Let $\tilde q$ be obtained from $q$ by dropping the first component.
Define $m_{11}$ as the element of $M$ in position $(1,1)$, $\tilde{m}$ as the first column of $M$ without the first element, and $\tilde{M}$ is the submatrix of $M$ obtained by stripping $M$ of its first row and first column:
$$
q=\left(\begin{array}{c}
q_1 \\
\hline \tilde{q}
\end{array}\right) ,
\quad
M=\left(\begin{array}{c|c}
m_{11} & \tilde{m}^\transp \\
\hline \tilde{m} & \tilde{M}
\end{array}\right) .
$$
We can now write $\E'$ as follows:
\begin{align*}
\E' &= \bra{\tilde y \in \R^{n-1} : (\tilde y - \tilde q)^\transp \tilde M (\tilde y - \tilde q) - 2 q_1 {\tilde m}^\transp (\tilde y - \tilde q) \le 1 - q_1^2 m_{11}} \\
& = \bra{\tilde y \in \R^{n-1} : \pare{\tilde y - c'}^\transp \tilde M \pare{\tilde y - c'} \le \zeta},
\end{align*}
where $c' :=\tilde q+q_1 \tilde{M}^{-1} \tilde{m}$ and $\zeta := 1-q_1^2 \pare{m_{11}-\tilde{m}^\transp \tilde{M}^{-1} \tilde{m}}$.
Observe that $\E' = \emptyset$ if $\zeta < 0$, $\E' = \bra{c'}$ if $\zeta = 0$, and $\E'$ is the ellipsoid $\E(c',\tilde M / \zeta)$ if $\zeta > 0$.
\end{prf}

We remark that, in \cref{lem mitr go to lower dim}, it is fundamental that ellipsoids are encoded by means of their symmetric positive definite matrix, i.e., the matrix $Q$ in $\E(c,Q)$.
A different encoding of the ellipsoid,
where one encodes a matrix $L$ such that $L^\transp L = Q$ instead of the matrix $Q$, would not allow us to obtain the ellipsoid $\E(c',Q')$ using only rational numbers.

We only need one more simple result, which illustrates a basic property of $\epsilon$-approximate solutions.
The proof is given for completeness.

\begin{lemma}
\label{lem partition}
Let $\S^1,\dots,\S^t \subseteq \R^n$.
Let $\genobj : \cup_{i=1}^t \S^i \to \R$ and assume that it has a (global) minimum.
Let $\epsilon \in [0,1]$.
For $i=1,\dots,t$, let $x^i$ be an $\epsilon$-approximate solution to the optimization problem $\inf\{\genobj(x) : x \in \S^i\}$.
Then, each optimal solution to $\min\{\genobj(x) : x \in \{x^1,\dots,x^t\}\}$, is an $\epsilon$-approximate solution to the optimization problem $\inf\{\genobj(x) : x \in \cup_{i=1}^t \S^i\}$.
\end{lemma}

\begin{prf}
Let $\genobj_{\inf}$ and $\genobj_{\sup}$ denote the infimum and the supremum of~$\genobj(x)$ on the domain $\cup_{i=1}^t \S^i$.
Let $x^*$ be an optimal solution to $\min\{\genobj(x) : x \in \cup_{i=1}^t \S^i\}$ and let $j \in \{1,\dots,t\}$ such that $x^* \in \S^j$.
Let $\genobj^j_{\inf}$ and $\genobj^j_{\sup}$ denote the infimum and the supremum of~$\genobj(x)$ on the set $\S^j$.
We then have
\begin{align*}
\genobj(x^j) - \genobj_{\inf} 
& = \genobj(x^j) - \genobj^j_{\inf} \\
& \le \epsilon (\genobj^j_{\sup} - \genobj^j_{\inf}) \\
& \le \epsilon (\genobj_{\sup} - \genobj_{\inf}).
\end{align*}
Here, in the first equality we used $\genobj_{\inf} = \genobj(x^*) = \genobj^j_{\inf}$, the first inequality holds because $x^j$ is an $\epsilon$-approximate solution to $\inf\{\genobj(x) : x \in \S^j\}$, and the last inequality holds because $\genobj^j_{\sup} \le \genobj_{\sup}$ and $\genobj^j_{\inf} = \genobj_{\inf}$.
Hence, $x^j$ is an $\epsilon$-approximate solution to $\inf\{\genobj(x) : x \in \cup_{i=1}^t \S^i\}$.
Therefore, also each optimal solution to $\min\{\genobj(x) : x \in \{x^1,\dots,x^t\}\}$ is an $\epsilon$-approximate solution to $\inf\{\genobj(x) : x \in \cup_{i=1}^t \S^i\}$.
\end{prf}


\subsection{Proof of \cref{th MITR algorithm}}
\label{sec MITR algorithm}

We are now ready to prove \cref{th MITR algorithm}.
Before presenting the proof, we give an overview of the algorithm.
First, we apply \cref{prop E-MIQP flat}.
If we find an $\epsilon$-approximate solution to \cref{prob E-MIQP}, we are done.
Otherwise, we find a vector $\direction \in \Z^n$ that allows us to partition the feasible region of \cref{prob E-MIQP} into a number of parallel hyperplanes that is polynomial in $1/\epsilon$.
Let $\H(\direction,\beta)$ be one of these hyperplanes and let \cref{prob E-MIQP beta} be obtained from \cref{prob E-MIQP} by adding the constraint $\direction^\transp x = \beta$.
If $\S:=\H(d,\beta) \cap ( \Z^p \times \R^{n-p} )$ or $\E(c,Q) \cap \H(d,\beta)$ is empty, then \cref{prob E-MIQP beta} is infeasible, so assume both sets are nonempty.
If $\E(c,Q) \cap \H(d,\beta)$ is a singleton, then we can easily solve \cref{prob E-MIQP beta}.
Otherwise, we use \cref{lem mitr go to lower dim} to transform \cref{prob E-MIQP beta} into an equivalent problem of the form of \cref{prob E-MIQP}, but with $p-1$ integer variables.
We can then apply the algorithm recursively.
We now give the complete proof of \cref{th MITR algorithm}, which uses \cref{prop E-MIQP flat,lem mitr go to lower dim,lem partition}.

\begin{prf}
We now present our approximation algorithm for \cref{prob E-MIQP}.

\noindent
\textbf{Step 1: Approximation or partition.}
We first apply the algorithm in \cref{prop E-MIQP flat}.
If the algorithm finds an $\epsilon$-approximate solution to \cref{prob E-MIQP}, then we are done.
Otherwise, the algorithm finds a nonzero vector $\direction \in \Z^n$ with $d_{p+1} = \cdots = d_n = 0$ and a scalar $\rho \in \Q$ such that 
\begin{align*}
\{d^\transp x : x \in \E(c,Q)\} \subseteq \sbra{\rho, \rho + 19 p \constlen / \epsilon}.
\end{align*}
Note that $p \in \{1,\dots,n\}$ and $\direction^\transp x$ is an integer for every $x \in \Z^p \times \R^{n-p}$.
Therefore, every feasible point of \cref{prob E-MIQP} is contained in one of the hyperplanes
\begin{align*}
\{x \in \R^n : {\direction}^\transp x = \beta\}, \qquad \beta = \ceil{\rho}, \ \ceil{\rho} + 1, \dots, \floor{\rho + 19 p \constlen / \epsilon}.
\end{align*}
For each $\beta = \ceil{\rho}, \ceil{\rho} + 1, \dots, \floor{\rho + 19 p \constlen / \epsilon}$, we define the optimization problem
\begin{align}
\label[problem]{prob E-MIQP beta}
\tag{E-MIQP-$\beta$}
\begin{split}
\min & \quad x^\transp H x + h^\transp x \\
\st & \quad (x-c)^\transp Q (x-c) \leq 1 \\
& \quad \direction^\transp x = \beta \\
& \quad x \in \Z^p \times \R^{n-p}.
\end{split}
\end{align}

If $n=1$, \cref{prob E-MIQP beta} can be trivially solved, thus we now assume $n \ge 2$.
For each $\beta$, we apply the algorithm in \cref{lem mitr go to lower dim}.
This algorithm determines whether the sets $\S:=\H(d,\beta) \cap ( \Z^p \times \R^{n-p} )$ and $\E(c,Q) \cap \H(d,\beta)$ are empty or not.
If any of the two is empty, then \cref{prob E-MIQP beta} is infeasible, thus we now assume that both sets are nonempty. 
Then, the algorithm finds 
a map $\eta: \R^{n-1} \to \R^n$ of the form $\eta(y) := \bar x + Ty$, with $\bar x \in \Z^n$ and $T \in \Q^{n\times {n-1}}$ of full rank, 
and the preimage $\E'$ of the set $\E(c,Q) \cap \H(\direction,\beta)$ under $\eta$.
Assume first that $\E'$ is a singleton $\bra{c'}$, for some $c' \in \Q^{n-1}$.
Then $\E(c,Q) \cap \H(\direction,\beta)$ is the singleton $\bra{\bar x + Tc'}$.
We check if $\bar x + Tc' \in \Z^p \times \R^{n-p}$.
If it is, then it is the optimal solution to \cref{prob E-MIQP beta}, and if it is not then \cref{prob E-MIQP beta} is infeasible.
Thus we now assume that $\E'$ is an ellipsoid $\E(c',Q') \subset \R^{n-1}$.
We then have 
\begin{align*}
\H(d,\beta) & =  \eta(\R^{n-1}), \\
\S & =
\eta(\Z^{p-1} \times \R^{n-p}), \\
\E(c,Q) \cap \H(d,\beta) & =  \eta( \E(c',Q')).
\end{align*}
Applying the change of variables $x = \bar x + Ty$ to \cref{prob E-MIQP beta}, we obtain the optimization problem
\begin{align}
\label[problem]{prob MITR beta'}
\tag{MITR-$\beta'$}
\begin{split}
\gamma
+ 
\min & \quad y^\transp H' y + {h'}^\transp y \\
\st & \quad y \in \E(c',Q') \\
& \quad x \in \Z^{p-1} \times \R^{n-p},
\end{split}
\end{align}
where $\gamma:= \bar x^\transp H \bar x + {h}^\transp \bar x$, $H' := T^{\transp} H T$, $h' := T^\transp {h} + 2 T^\transp H^\transp \bar x$.
Note that \cref{prob MITR beta'} has the form of \cref{prob E-MIQP}, but with $p-1$ integer variables and $n-p$ continuous variables.

\noindent
\textbf{Step 2: Recursion.}
We apply recursively the algorithm to each obtained \cref{prob MITR beta'}, for each $\beta = \ceil{\rho}, \ceil{\rho} + 1, \dots, \floor{\rho + 19 p \constlen / \epsilon}$.
Since each time the problem is partitioned in Step~1, the number of integer variables decreases by one, the total number of iterations of the algorithm is upper bounded by 
$$
\pare{19 p \constlen / \epsilon}^p .
$$

\noindent
\textbf{Step 3: Termination.}
If the algorithm did not find any feasible solution to \cref{prob E-MIQP}, the algorithm correctly returns that the feasible region is empty.
Otherwise, among all the feasible solutions to \cref{prob E-MIQP} found, the algorithm returns one with the minimum objective value.
The result then follows from \cref{lem partition}.
This is because the function $f(x) := x^\transp H x + h^\transp x$, from $\E(c,Q) \cap (\Z^p \times \R^{n-p})$ to $\R$, is a continuous function from a compact set to the real numbers.
Hence, from the extreme value theorem, it has a (global) minimum.
\end{prf}

\section{Mixed integer quadratic programming}
\label{sec MIQP}





In this section we prove \cref{th MIQP algorithm}.

\subsection{Flatness result for \cref{prob MIQP}}

Our first goal is to present a flatness result for \cref{prob MIQP}:
an algorithm that either finds an approximate solution, or finds a nonzero vector $d$ along which the feasible region is flat.
This result has several similarities with the flatness result for \cref{prob E-MIQP} that we presented in \cref{prop E-MIQP flat}.
The key difference, besides the optimization problem considered, is that the algorithm in \cref{prop E-MIQP flat} can be applied with any $\epsilon \in (0,1]$, while the algorithm in the flatness result for \cref{prob MIQP} only works with $\epsilon = 1-\Theta(n^{-2})$.


\begin{proposition}[Flatness result for \cref{prob MIQP}]
\label{prop MIQP flat}
Consider \cref{prob MIQP}, and assume that $\{x \in \R^n : Wx \le w\}$ is full-dimensional and bounded.
There is an algorithm which either finds an $\epsilon$-approximate solution,
where the approximation factor is $\epsilon = 1-\Theta(n^{-2})$, 
or finds a nonzero vector $\direction \in \Z^n$ with $d_{p+1} = \cdots = d_n = 0$ and a scalar $\rho \in \Q$ such that 
\begin{align*}
\bra{d^\transp x : Wx \le w} \subseteq \sbra{\rho, \rho + 14 (n+1) p \constlen}.
\end{align*}
The running time of the algorithm is polynomial in the size of \cref{prob MIQP}.
\end{proposition}


Before giving the proof, we give an overview of the algorithm.
The structure is similar to the algorithm in the proof of \cref{prop E-MIQP flat}, with few key differences.
First, we construct an affine transformation $\tau$ that maps $\P$ into a polytope that is sandwiched between $\B(0,1)$ and $\B(0,n+1)$.
We denote by \cref{prob MIQP'} the optimization problem obtained by applying the same affine transformation to \cref{prob MIQP}.
Next, we consider \cref{prob MIQP'-S} obtained from \cref{prob MIQP'} by replacing the constraint $y \in \tau( \P )$ with with the stronger requirement $y \in \B(0,1)$.
Using \cref{prop MITR dist points}, we search for two feasible solutions $y^\#,y^\diamond$ to \cref{prob MIQP'-S} such that $f(y^\#)$ is sufficiently larger than $y^\diamond$.
If \cref{prop MITR dist points} finds the vectors $y^\#,y^\diamond$, we show that $y^\diamond$ is an $\epsilon$-approximate solution to \cref{prob MIQP'}, where the approximation factor is $\epsilon = 1-\Theta(n^{-2})$.
Otherwise, we show that $\P$ is flat and we find the associated vector $d$.
We now give the complete proof of \cref{prop MIQP flat}, which uses \cref{lem sandwich ellipsoid balls,prop MITR dist points,lem direction transformation}.

\begin{prf}
Let $\P := \{x \in \R^n : Wx \le w\}$.
Using theorem~15.6 in~\cite{SchBookIP} (with $\gamma:=1/(n+1/2)$), we find a rational ellipsoid $\E(c,Q)$ such that 
\begin{align*}
\E(c,(n+1/2)^2 Q) \subseteq \P \subseteq \E(c,Q).
\end{align*}
Using the algorithm in \cref{lem sandwich ellipsoid balls} (with $Q := (n+1/2)^2 Q$ and $\delta := 1/(2n+1)$), we find a map $\tau: \R^n \to \R^n$ of the form $\tau(x) := B(x-c)$, with $B \in \Q^{n\times n}$ invertible, such that 
\begin{align*}
\B(0,1) \subseteq \tau( \E(c,(n+1/2)^2 Q) ) \subseteq \B(0,(2n+2)/(2n+1)).
\end{align*}
In particular, the first containment implies $\B(0,1) \subseteq \tau( \P )$.
Next, we prove $\tau( \P ) \subseteq \B(0,n+1)$:
\begin{align*}
\tau( \P )
& \subseteq \tau( \E(c,Q) ) \\
& = \E(0,B^{-\transp} Q B^{-1}) \\
& = (n+1/2) \ \E(0,(n+1/2)^2 B^{-\transp} Q B^{-1}) \\
& = (n+1/2) \ \tau( \E(c,(n+1/2)^2 Q) ) \\
& \subseteq (n+1/2) \ \B(0,(2n+2)/(2n+1)) \\
& = \B(0,n+1).
\end{align*}
Therefore, we have shown
\begin{align*}
\B(0,1) \subseteq \tau( \P ) \subseteq \B(0,n+1).
\end{align*}
%
%
Applying the change of variables $y := B (x-c)$ to \cref{prob MIQP}, we obtain the optimization problem
\begin{align}
\label[problem]{prob MIQP'}
\tag{MIQP$'$}
\begin{split}
\gamma
+ 
\min & \quad y^\transp H' y 
+ {h'}^\transp y
 \\
\st & \quad y \in \tau( \P ) \\
& \quad y \in \Pi_p(b^1,\dots,b^n) + \{c'\},
\end{split}
\end{align}
where $\gamma:= c^\transp H c + {h}^\transp c$, $H' := B^{-\transp} H B^{-1}$, $h' := 2 B^{-\transp} H^\transp c + B^{-\transp} h$, $b^1,\dots,b^n$ are the columns of $B$, and $c':=-Bc$.
For ease of notation, in the following we denote by $f_y(y) := y^\transp H' y + {h'}^\transp y$.



Consider now the mixed integer trust region problem obtained from \cref{prob MIQP'} by replacing the constraint $y \in \tau( \P )$ with the stronger requirement $y \in \B(0,1)$:
\begin{align}
\label[problem]{prob MIQP'-S}
\tag{MIQP$'$-s}
\begin{split}
\gamma + 
\min & \quad y^\transp H' y 
+ {h'}^\transp y
 \\
\st & \quad y^\transp y \le 1 \\
& \quad y \in \Pi_p(b^1,\dots,b^n) + \{c'\}.
\end{split}
\end{align} 
Using the algorithm in \cref{prop MITR dist points} (with $\epsilon := 1/14$), we either find two feasible solutions $y^\#,y^\diamond$ to \cref{prob MIQP'-S} such that
\begin{align}
\label{eq MIQP from prop MITR dist points}
f(y^\#)-f(y^\diamond) \geq 
\frac 12 \max \bra{\norm{H'},\norm{h'}},
\end{align}
or we find a nonzero vector 
${\direction'} \in \Q^n$ with ${\direction'}^\transp b^1, \dots, {\direction'}^\transp b^p$ integer and ${\direction'}^\transp b^{p+1} = \cdots = {\direction'}^\transp b^n = 0$ such that $\width_{\direction'}(\B(0,1)) \le 14 p \constlen$.

In the remainder of the proof, we consider separately two cases.
First, consider the case where \cref{prop MITR dist points} successfully found the two vectors $y^\#,y^\diamond$.
Note that $y^\#$ and $y^\diamond$ are feasible also to \cref{prob MIQP'}.
Denote by $f_{\sup}$ and $f_{\inf}$ the supremum and the infimum of $f_y(y)$ on the feasible region of \cref{prob MIQP'}.
From \eqref{eq MIQP from prop MITR dist points}, we directly obtain the following lower bound on $f_{\sup}-f(y^\diamond)$:
\begin{align}
\label{eq MIQP lower bound}
\begin{split}
f_{\sup}-f(y^\diamond) & \ge f(y^\#)-f(y^\diamond) \\
& \ge \frac 12 \max \bra{\norm{H'},\norm{h'}}.
\end{split}
\end{align}
%
Next, we obtain the following upper bound on $f_{\sup}-f_{\inf}$:
\begin{align}
\label{eq MIQP upper bound}
f_{\sup}-f_{\inf} \leq 4 (n+1)^2 \max \bra{\norm{H'}, \norm{h'}}.
\end{align}
To prove this bound, observe that, for any $y \in \B(0,n+1)$,
\begin{align*}
|f(y)| & \le \abs{y^\transp H' y} + \abs{{h'}^\transp y} \\
& \le \norm{H'} \norm{y}^2 + \norm{h'} \norm{y} \\
& \leq (n+1)^2 \norm{H'}+ (n+1) \norm{h'} \\
&\le 2 (n+1)^2 \max \bra{\norm{H'},\norm{h'}}.
\end{align*}

We are now ready to show that $y^\diamond$ is an $\epsilon$-approximate solution to \cref{prob MIQP'}, where the approximation factor is $\epsilon = 1-\Theta(n^{-2})$.
We combine \eqref{eq MIQP lower bound} and \eqref{eq MIQP upper bound} to obtain:
\begin{align*}
f_{\sup}-f_{\inf}
\le
4 (n+1)^2 \max \bra{\norm{H'}, \norm{h'}}
\le 
8 (n+1)^2 \pare{f_{\sup}-f(y^{\diamond})}.
\end{align*}
We can then use the above inequality to obtain the bound:
\begin{align*}
f(y^{\diamond}) - f_{\inf} 
& = (f_{\sup} - f_{\inf}) -  (f_{\sup} - f(y^{\diamond})) \\
& \le
\pare{1-\frac{1}{8(n+1)^2}}
\pare{f_{\sup}-f_{\inf}}.
\end{align*}
It then follows that $x^\diamond := B^{-1} y^\diamond + c$ is an $\epsilon$-approximate solution to \cref{prob MIQP}, where the approximation factor is $\epsilon = 1-\Theta(n^{-2})$.
This concludes the first case.

Next, consider the case where \cref{prop MITR dist points} found a vector 
${\direction'} \in \Q^n$ with ${\direction'}^\transp b^1, \dots, {\direction'}^\transp b^p$ integer and ${\direction'}^\transp b^{p+1} = \cdots = {\direction'}^\transp b^n = 0$ such that
$
\width_{\direction'}(\B(0,1)) \le 14 p \constlen.
$
%
We can now upper bound $\width_{\direction'}\pare{\B(0,n+1)}$ as follows:
\begin{align*}
\width_{\direction'}\pare{\B(0,n+1)} 
& = (n+1) \width_{\direction'}\pare{\B(0,1)} \\
& \le 14 (n+1) p \constlen.
\end{align*}


Now let $\direction := B^\transp {\direction'} \in \Q^n$.
From \cref{lem direction transformation}, $d \in \Z^n$, $d_{p+1} = \cdots = d_n = 0$, and, for any nonempty bounded closed set $\S \subseteq \R^n$, $\width_\direction(\S) = \width_{\direction'}(\tau(\S)).$
If we let $\S := \tau^\leftarrow(\B(0,n+1))$, where $\tau^\leftarrow$ denotes the inverse of $\tau$, we obtain
\begin{align*}
\width_\direction(\tau^\leftarrow(\B(0,n+1))) 
= \width_{\direction'}(\B(0,n+1))
\le 14 (n+1) p \constlen.
\end{align*}
%
Let $\rho := \min\{{\direction}^\transp x : x \in \tau^\leftarrow(\B(0,n+1))\}$, and note that we can calculate it as follows:
\begin{align*}
\rho 
& = \min\{{\direction}^\transp (B^{-1} y + c) : y \in \B(0,n+1)\} \\
& = \min\{{\direction'}^\transp y : y \in \B(0,n+1)\} + \direction^\transp c \\
& = - (n+1) \norm{\direction'} + \direction^\transp c.
\end{align*}
We obtain
\begin{align*}
\{d^\transp x : x \in \tau^\leftarrow(\B(0,n+1))\} \subseteq [\rho, \rho + 14 (n+1) p \constlen].
\end{align*}
Since $\P \subseteq \tau^\leftarrow(\B(0,n+1))$, we also have
\begin{align*}
\{d^\transp x : x \in \P\} \subseteq [\rho, \rho + 14 (n+1) p \constlen].
\end{align*}
This concludes the second case.
\end{prf}

\subsection{Approximation algorithm for \cref{prob MIQP}}

We are now ready to prove \cref{th MIQP algorithm}.
Before presenting the proof, we give an overview of the algorithm. 
First, we reduce ourselves to the case where the polyhedron $\bra{x \in \R^n : Wx \le w}$ is bounded and full-dimensional.
Next, we apply \cref{prop MIQP flat}.
If we find an $\epsilon$-approximate solution to \cref{prob MIQP} where the approximation factor is $\epsilon = 1-\Theta(n^{-2})$, then we are done.
Otherwise, we find a vector $\direction \in \Z^n$ that allows us to partition the feasible region of \cref{prob MIQP} into a number of parallel hyperplanes that is polynomial in $n$.
Let $\H(\direction,\beta)$ be one of these hyperplanes and let \cref{prob MIQP beta} be obtained from \cref{prob MIQP} by adding the constraint $\direction^\transp x = \beta$.
We transform \cref{prob MIQP beta} into an equivalent problem of the form of \cref{prob MIQP}, but with $p-1$ integer variables.
We can then apply the algorithm recursively.
We now give the complete proof of \cref{th MIQP algorithm}, which uses \cref{prop MIQP flat,lem partition}.


\begin{prf}
We now present our approximation algorithm for \cref{prob MIQP}.

\noindent
\textbf{Step 0: Bounded feasible region.}
Consider \cref{prob MIQP} and assume there exists $\bar f \in \R$ such that every feasible solution has objective value at least $\bar f$.
Theorem~4 in~\cite{dPDeyMol17MPA} implies that, if the feasible region in nonempty, there is an optimal solution of size bounded by an integer $\psi$, which is polynomial in the size of the problem. \footnote{Even though Theorem~4 in~\cite{dPDeyMol17MPA} does not give $\psi$ explicitly, a formula for $\psi$, as a function of the size of \cref{prob MIQP}, can be derived from its proof.}
We add to the system $Wx \le w$ the linear inequalities $-2^\psi \le x_i \le 2^\psi$, for $i=1,\dots, n$, whose size is polynomial in the size of \cref{prob MIQP}.
It is simple to check that an $\epsilon$-approximate solution to the obtained problem is also an $\epsilon$-approximate solution to the original, for every $\epsilon \in [0,1]$.
Therefore, we can now assume that $\bra{x \in \R^n : Wx \le w}$ is bounded.


\noindent
\textbf{Step 1: Full dimensionality.}
Consider the sets 
\begin{align*}
\P := \bra{x \in \R^n : Wx \le w}, \qquad \S := \P \cap \pare{\Z^p \times \R^{n-p}},
\end{align*}
and note that $\S$ constitutes the feasible region of \cref{prob MIQP}.
We apply the algorithm in theorem~1 in \cite{dP24mSIOPT}.
The algorithm either returns that $\S$ is empty, or it returns $p' \in \bra{0,1,\dots,p}$, $n' \in \bra{p',p'+1,\dots,p'+n-p}$, 
a map $\tau : \R^{n'} \to \R^n$
of the form $\tau(x') = \bar x + Mx'$, with 
$\bar x \in \Z^p \times \Q^{n-p}$ and $M \in \Q^{n \times n'}$ of full rank, such that
the polyhedron 
\begin{align*}
\P' := \bra{x' \in \R^{n'} : WMx' \le w - W \bar x}
\end{align*}
is full-dimensional, and 
\begin{align*}
\P & = \tau \pare{\P'} \\
\S & = \tau \pare{\P' \cap \pare{\Z^{p'} \times \R^{n'-p'}}}.
\end{align*}
%
If $\S$ is empty we are done.
Otherwise, we perform the change of variables $x = \bar x + M x'$ to \cref{prob MIQP} and, after dropping the constant $\bar x^\transp H \bar x + h^\transp \bar x$ in the objective function, we obtain the problem
\begin{align}
\label[problem]{prob MIQP fulldim}
\tag{MIQP$'$}
\begin{split}
\min & \quad {x'}^\transp H' x' + {h'}^\transp x' \\
\st & \quad W'x' \le w' \\
& \quad x' \in \Z^{p'} \times \R^{n'-p'},
\end{split}
\end{align} 
where $H' := M^\transp H M$,
${h'}^\transp := h^\transp M + 2 \bar x^\transp HM$,
$W' := WM$, 
and $w' := w - W \bar x$.
Since $\P$ is bounded, also the polyhedron $\P'$, which is full-dimensional, is also bounded.

Since the definition of $\epsilon$-approximate solution is preserved under changes of variables and translations of the objective function, it then suffices to find an $\epsilon$-approximate solution to \cref{prob MIQP fulldim}.
%
%
For ease of notation, 
from now on we consider \cref{prob MIQP} instead of \cref{prob MIQP fulldim} and we further assume, without loss of generality, that $\{x \in \R^n : Wx \le w\}$ in \cref{prob MIQP} is full-dimensional.

\noindent
\textbf{Step 2: Approximation or partition.}
We now apply the algorithm in \cref{prop MIQP flat}.
If the algorithm finds an $\epsilon$-approximate solution to \cref{prob MIQP} where the approximation factor is $\epsilon = 1-\Theta(n^{-2})$, then we are done.
Otherwise, the algorithm finds a nonzero vector $\direction \in \Z^n$ with $d_{p+1} = \cdots = d_n = 0$ and a scalar $\rho \in \Q$ such that 
\begin{align*}
\bra{d^\transp x : Wx \le w} \subseteq \sbra{\rho, \rho + 14 (n+1) p \constlen}.
\end{align*}
Note that $p \in \{1,\dots,n\}$ and $\direction^\transp x$ is an integer for every $x \in \Z^p \times \R^{n-p}$.
Therefore, every feasible point of \cref{prob MIQP} is contained in one of the hyperplanes
\begin{align*}
\{x \in \R^n : {\direction}^\transp x = \beta\}, \qquad \beta = \ceil{\rho}, \ \ceil{\rho} + 1, \dots, \floor{\rho + 14 (n+1) p \constlen}.
\end{align*}
For each $\beta = \ceil{\rho}, \ceil{\rho} + 1, \dots, \floor{\rho + 14 (n+1) p \constlen}$, we define the optimization problem
\begin{align}
\label[problem]{prob MIQP beta}
\tag{MIQP-$\beta$}
\begin{split}
\min & \quad x^\transp H x + h^\transp x \\
\st & \quad Wx \le w \\
& \quad \direction^\transp x = \beta \\
& \quad x \in \Z^p \times \R^{n-p}.
\end{split}
\end{align}
If $n=1$, \cref{prob MIQP beta} can be trivially solved, thus we now assume $n \ge 2$.
For each $\beta$, we apply theorem~1 in \cite{dP24mSIOPT} to the feasible region of \cref{prob MIQP beta}.
If we find out that the feasible region is empty, we are done.
Otherwise, we transform, with a change of variables, \cref{prob MIQP beta} into a new problem, which we denote by Problem~MIQP-$\beta'$.
Since $d_{p+1} = \cdots = d_n = 0$, Problem~MIQP-$\beta'$ has at most $p-1$ integer variables and $n-p$ continuous variables.


\noindent
\textbf{Step 3: Recursion.}
We apply recursively the algorithm, from Step~1, to each obtained Problem~MIQP-$\beta'$, for each $\beta = \ceil{\rho}, \ceil{\rho} + 1, \dots, \floor{\rho + 14 (n+1) p \constlen}$.
Since each time the problem is partitioned in Step~2, the number of integer variables decreases by one, the total number of iterations of the algorithm is upper bounded by 
$$
\pare{14 (n+1) p \constlen}^p.
$$

\noindent
\textbf{Step 4: Termination.}
If the algorithm did not find any feasible solution to \cref{prob MIQP}, the algorithm correctly returns that the feasible region is empty.
Otherwise, among all the feasible solutions to \cref{prob MIQP} found, the algorithm returns one with the minimum objective value.
The result then follows from \cref{lem partition}.
This is because the function $f(x) := x^\transp H x + h^\transp x$, from the original feasible region $\{x \in \Z^p \times \R^{n-p} : Wx \le w\}$ to $\R$, is lower bounded by assumption, and thus has a (global) minimum, due to theorem~4 in~\cite{dPDeyMol17MPA}.
\end{prf}



\bigskip

\begin{small}
\noindent
\textbf{Funding:} A.~Del~Pia is partially funded by AFOSR grant FA9550-23-1-0433. Any opinions, findings, and conclusions or recommendations expressed in this material are those of the authors and do not necessarily reflect the views of the Air Force Office of Scientific Research.

\bigskip

\noindent
\textbf{Acknowledgments:} The author would like to thank two anonymous referees for comments and suggestions that improved the quality of this manuscript.
\end{small}


\ifthenelse {\boolean{MPA}}
{
\bibliographystyle{spmpsci}
}
{
\bibliographystyle{plain}
}


\end{document}